\documentclass[%
reprint,
superscriptaddress,https://www.overleaf.com/project/60418f9384090c596f682ce0
 nofootinbib,
 amsmath,amssymb,
 aps,
 pre,
]{revtex4-2}

\usepackage{graphicx}% Include figure files
\usepackage{dcolumn}% Align table columns on decimal point
\usepackage{bm}% bold math
\usepackage{bbm}% bold math
\usepackage{float}
\usepackage{subcaption}
\usepackage{multirow}
\usepackage{comment}
\usepackage{dsfont}
\usepackage{amsthm}
\usepackage{amsmath, esint}
\makeatletter
\renewcommand*\env@matrix[1][*\c@MaxMatrixCols c]{%
  \hskip -\arraycolsep
  \let\@ifnextchar\new@ifnextchar
  \array{#1}}
\makeatother
\usepackage{physics}
\usepackage{xcolor}
\usepackage{makecell}
\usepackage{scalerel,stackengine}
\usepackage{booktabs}
\usepackage{kantlipsum}
\allowdisplaybreaks

\newcommand{\mat}[1]{\mathbf{#1}}

\def\spose#1{\hbox to 0pt{#1\hss}}
\def\simlt{\mathrel{\spose{\lower 3pt\hbox{$\mathchar"218$}}
     \raise 2.0pt\hbox{$\mathchar"13C$}}}
\def\simgt{\mathrel{\spose{\lower 3pt\hbox{$\mathchar"218$}}
     \raise 2.0pt\hbox{$\mathchar"13E$}}}
\def\simpropto{\mathrel{\spose{\lower 3pt\hbox{$\mathchar"218$}}
     \raise 2.0pt\hbox{$\propto$}}}

\def\beq#1{\begin{equation}\label{#1}}
\def\eeq{\end{equation}}
\def\beqa#1{\begin{eqnarray}\label{#1}}
\def\eeqa{\end{eqnarray}}

\begin{document}

%\preprint{APS/123-QED}

%\title{Discovering Conservation Laws Only Requires Linear Algorithms}

\title{Discovering New Interpretable Conservation Laws as Sparse Invariants}

%\author{Ziming Liu, Obin Sturm, Saketh Bharadwaj, Sam Silva, Max Tegmark}
\author{Ziming Liu}
\affiliation{
Department of Physics, Institute of Artificial Intelligence and Fundamental Interactions, Massachusetts Institute of Technology, Cambridge, USA}
\author{Patrick Obin Sturm}
\affiliation{ Department of Earth Sciences, University of Southern California, Los Angeles, USA
% obin, please fill your institution here
}
\author{Saketh Bharadwaj}
\affiliation{Department of Chemical Engineering, Indian Institute of Technology, Hyderabad, India
% Saketh, please fill your institution here.
}
\author{Sam J. Silva}
\affiliation{ Department of Earth Sciences, University of Southern California, Los Angeles, USA
% Obin or Sam, please fill Sam's instituion here (the same as Obin's?)
}
\author{Max Tegmark}
\affiliation{
Department of Physics, Institute of Artificial Intelligence and Fundamental Interactions, Massachusetts Institute of Technology, Cambridge, USA}

\date{\today}% It is always \today, today,
             %  but any date may be explicitly specified

\begin{abstract}
Discovering conservation laws for a given dynamical system is important but challenging.  In a \textit{theorist} setup (differential equations and basis functions are both known), we propose the {\bf S}parse {\bf I}nvariant {\bf D}etector (SID), an algorithm that auto-discovers conservation laws from differential equations. Its algorithmic simplicity allows robustness and interpretability of the discovered conserved quantities. We show that SID is able to rediscover known and even discover new conservation laws in a variety of systems. For two examples in fluid mechanics and atmospheric chemistry, SID discovers 14 and 3 conserved quantities, respectively, where only 12 and 2 were previously known to domain experts.
\end{abstract}

%\keywords{Suggested keywords}%Use showkeys class option if keyword
                              %display desired
\maketitle

%\onecolumngrid

\section{Introduction}

Conservation laws are important concepts in physics, yet discovering them is challenging. Ideally, the set of discovered conserved quantities should be \textit{complete}, \textit{independent} and \textit{interpretable}. Although several attempts have been made to automate the discovery process with machine learning ~\cite{poincare1,mototake2021interpretable,wetzel2020discovering,PhysRevResearch.3.L042035,arora2023model,kasim2022constants,poincare2,kaiser2018discovering}, their complicated setups and blackbox nature make it hard to guarantee all these desirable properties. This paper considers a simple yet realistic setup where all these desirable properties can be met. 

``Discovering conservation laws" can mean wildly different things for \textit{experimentalists}, \textit{computationalists} and \textit{theorists}, as shown in Table \ref{tab:setups}. Most prior work ~\cite{poincare1,mototake2021interpretable,wetzel2020discovering,PhysRevResearch.3.L042035,arora2023model,
  kasim2022constants}
takes on the experimentalist setup, assuming knowledge of neither the differential equations nor the form of conservation laws. ~\cite{poincare2} takes the computationalist setup,  assuming knowledge of differential equations.  This work explores the theorist setup, where both differential equations and basis functions of conservation laws are known. Admittedly, this setup is simpler than the other two, but is still realistic when theorists have the differential equations at hand and have educated guesses about the basis functions that may span the conserved quantities.

We propose {\bf S}parse {\bf I}nvariant {\bf D}etector (SID), an algorithm that reveals conservation Laws. SID is incredibly simple in the sense that it only requires linear algorithms (except for sparsification), so the results are much more trustworthy and interpretable than blackbox machine learning methods. 
%The closest prior work to ours is~\cite{poincare2}. However, they use neural networks to parameterize conserved quantities, so the reliability and interpretability are not guaranteed. 
Note that SID does not replace us human scientists, but rather acts as a helpful assistant: while humans need to input basis functions (i.e., \textit{formulating} hypotheses) to SID, SID is good at computing conserved quantities (i.e., \textit{testing} hypotheses) based on the given prompt. In this manner, human scientists can focus on the more creative part of the job, while SID does the technical and tedious work. %If necessary, multiple rounds of interactions between humans and SID are needed to jointly discover new conserved quantities.
This paper gives two examples where new conserved quantities are successfully discovered by SID: one in fluid mechanics, and another in atmospheric chemistry (see Table~\ref{tab:num_new_cq}). In the former one, although the new conserved quantities are somewhat expected in hindsight, humans alone may need several more months to find them. In the latter one, a new conserved quantity is found, which was unintended in the design of the model.  %If necessary, multiple rounds of interactions between humans and SID may be needed to jointly discover new conserved quantities.

% which was unintended in the design of the atmospheric chemistry model 

%Moreover, although all the prior works~\cite{poincare1,mototake2021interpretable,wetzel2020discovering,PhysRevResearch.3.L042035,arora2023model,kasim2022constants,poincare2} focus on benchmarking their algorithms by rediscovering known conservation laws, our method surprisingly discovers new conservation laws for two systems: atmospheric chemistry and fluid mechanics (3D), displayed in Table~\ref{tab:num_new_cq}

\begin{table}[tbp]
    \centering
    \caption{Three setups of conservation law discovery}
    \begin{tabular}{cccc}
    \toprule
     Setup  & Experimentalist & Computationalist & Theorist \\
     \midrule
    Model-based & No & Yes & Yes \\
    Known basis & No & No & Yes \\
    \midrule
    Independence & Partial & Yes & Yes \\
    Completeness & No & Partial & Yes \\
    Interpretability & Partial & Partial & Yes \\
    \midrule
    Reference & ~\cite{poincare1,mototake2021interpretable,wetzel2020discovering,PhysRevResearch.3.L042035,arora2023model,
  kasim2022constants} & ~\cite{poincare2} & This work \\
    \bottomrule
    \end{tabular}
    \label{tab:setups}
\end{table}

\section{Method}

\begin{figure*}[htbp]
    \centering
    \includegraphics[width=0.85\linewidth]{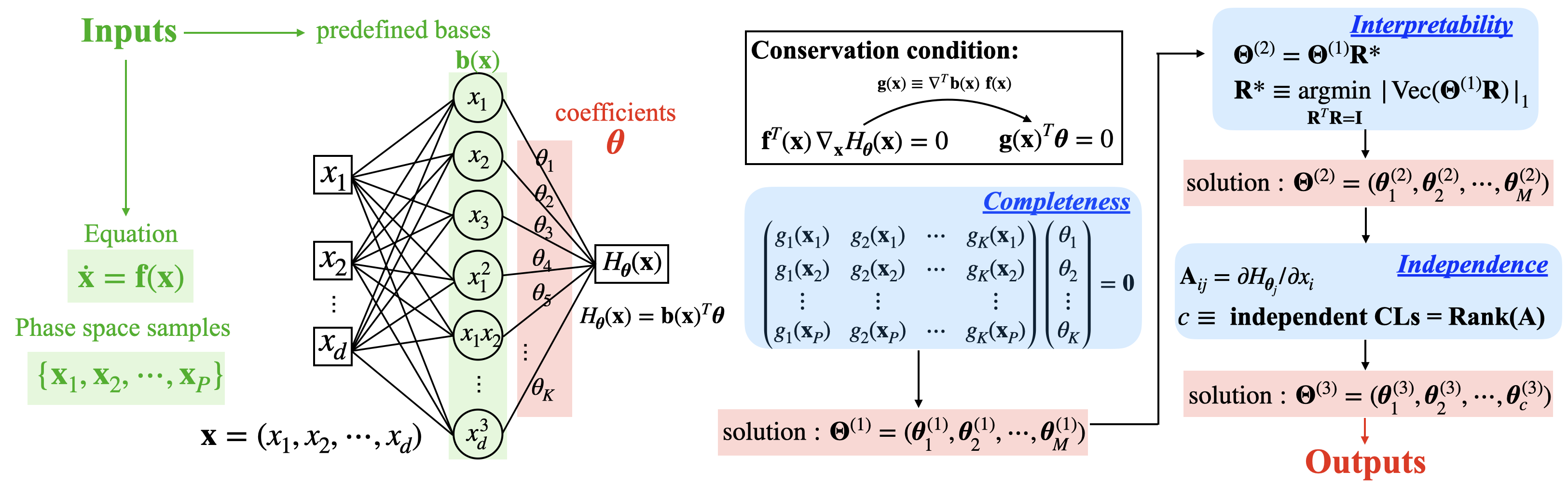}
    \caption{SID workflow: {\bf Inputs} are differential equations, basis functions and sample points; {\bf Outputs} are a set of conserved quantities which are complete, independent and interpretable.}
    \label{fig:algorithm}
\end{figure*}

\begin{table}[tbp]
    \centering
    \caption{The number of conserved quantities known to experts and discovered by SID}
    \begin{tabular}{|c|c|c|c|}\hline
         & Fluid (2D) & Fluid (3D)  & Atmosphere \\\hline
     Known  & {\bf 8} & 12 & 2 \\\hline 
     SID & {\bf 8} (simpler) & {\bf 14} & {\bf 3}  \\\hline
    \end{tabular}
    \label{tab:num_new_cq}
\end{table}

\subsection{Problem setup}
We consider a first-order differential equation $\dot{\mat{x}}=\mat{f}(\mat{x})$, where $\dot{\mat{x}}\equiv\frac{d\mat{x}}{dt}$, $\mat{x}\equiv(x_1,\cdots,x_d)\in\mathbb{R}^d$ is the state vector and $\mat{f}:\mathbb{R}^d\to\mathbb{R}^d$ is a vector field. This ODE formulation is more general than it seems: (1) Hamiltonian systems are subsumed as $\mat{x}\equiv(\mat{x}',\mat{p}')$; (2) Higher-order differential equations (e.g., $\ddot{\mat{y}}=\mat{f}(\mat{y})$) are included as $\mat{x}\equiv (\mat{y},\dot{\mat{y}},\cdots)$; (3) Partial differential equations (PDEs) become ODEs once discretized. 

A \textit{conserved quantity} is a scalar function $H(\mat{x}):\mathbb{R}^d\to\mathbb{R}$, such that its value remains constant along any trajectory obeying $\dot{\mat{x}}=\mat{f}(\mat{x})$~\footnote{The constant itself can be different for different trajectories.}. As proved in~\cite{poincare2}, a
necessary and sufficient condition for
$H(\mat{x})$ being a conserved quantity is $\nabla H(\mat{x})\cdot \mat{f}(\mat{x})=0$, since
\begin{equation}\label{eq:conservation-condition}
    0=\dot{H} = \nabla H(\mat{x})\cdot\dot{\mat{x}} = \nabla H(\mat{x})\cdot \mat{f}(\mat{x}).
\end{equation}

Given the differential equation $\dot{\mat{x}}=\mat{f}(\mat{x})$, we hope to find a set of conserved quantities $\{H_1,\cdots, H_c\}$ which satisfies these three properties:
\begin{itemize}
    \item \textit{Independence}: they are functionally independent, i.e., $g(H_1,\cdots,H_c)=0\Rightarrow g=0$.
    \item \textit{Completeness}: any conserved quantity $H$ (in the function space spanned by basis functions) can be expressed by them, i.e., there exists $g$ such that $H=g(H_1,\cdots,H_c)$.
    \item \textit{Interpretability}: conserved quantities can be written as (hopefully simple) symbolic formulas.
\end{itemize}

\subsection{Solving the linear equation and completeness}
The prior work~\cite{poincare2} parametrizes the conserved quantities  $H_{\boldsymbol\theta}(\mat{x})$ as neural networks and learns the parameters $\boldsymbol\theta$ to make $|\nabla H_{\boldsymbol\theta}(\mat{x})\cdot\mat{f}(\mat{x})|^2$ close to zero. However, neural network training may get stuck at local minima, so the results are not reliable. Moreover, the parameterized conserved quantities are not immediately interpretable. 

We consider a simpler setup. Assume that we know  $H_{\boldsymbol\theta}(\mat{x})$ to be a linear combination of $K$ predefined basis functions $b_i(\mat{x})\ (1\leq i\leq K)$ such that
\begin{equation}\label{eq:linear-regression}
    H_{\boldsymbol\theta}(\mat{x}) = \sum_{i=1}^K \theta_ib_i(\mat{x})\equiv{\boldsymbol\theta}\cdot \mat{b}(\mat{x}),
\end{equation}
where only $\boldsymbol\theta\in\mathbb{R}^K$ are learnable parameters to be determined and the vector $\mat{b}:\mathbb{R}^d\to\mathbb{R}^K$ defines the basis functions. Since the number of conserved quantities can exceed one, we define a set of parameters $\boldsymbol\Theta\equiv\{\boldsymbol\theta_1,\boldsymbol\theta_2,\cdots\}$ and their corresponding functions $H_{\boldsymbol\Theta}\equiv\{H_{\boldsymbol\theta}|\boldsymbol\theta\in\boldsymbol\Theta\}$. %A nice property of Eq.~(\ref{eq:linear-regression}) is that the function space and the parameter space have the same linear structure, i.e.,
%\begin{equation}\label{eq:linear-relation}
%\sum_i\alpha_iH_{\boldsymbol\theta_i}(\mat{x}) = H_{\sum_i \alpha_i{\boldsymbol\theta_i}}(\mat{x}).
%\end{equation}
As shown in FIG.~\ref{fig:algorithm}, Eq.~(\ref{eq:linear-regression}) is equivalent to a neural network whose last linear layer contains the only trainable parameters. With Eq.~(\ref{eq:linear-regression}), the conservation condition Eq.~(\ref{eq:conservation-condition}) becomes:
\begin{equation}\label{eq:simplified-condition}
    \mat{g}(\mat{x})^T{\boldsymbol\theta}=0,\quad \mat{g}(\mat{x})\equiv (\nabla \mat{b}(\mat{x}))\mat{f}(\mat{x}),
\end{equation}
which is a linear equation of ${\boldsymbol\theta}$. Remember that in our setup, both $\mat{b}(\mat{x})$ and $\mat{f}(\mat{x})$ are known, so $\mat{g}(\mat{x})\equiv(\nabla\mat{b}(\mat{x}))\mat{f}(\mat{x})$ is known as well. In practice, we draw $P$ random points $\mat{x}_i\ (1\leq i\leq P)$ from  phase space. A solution ${\boldsymbol\theta}$  should make Eq.~(\ref{eq:simplified-condition}) hold for all $\mat{x}_i$, or more explicitly,
\begin{equation}\label{eq:G_theta}
\underbrace{
    \begin{pmatrix}
    g_1(\mat{x}_1) & g_2(\mat{x}_1) & \cdots & g_K(\mat{x}_1) \\
    g_1(\mat{x}_2) & g_2(\mat{x}_2) & \cdots & g_K(\mat{x}_2) \\
    \vdots & \vdots & & \vdots \\
    g_1(\mat{x}_P) & g_2(\mat{x}_P) & \cdots & g_K(\mat{x}_P) \\
    \end{pmatrix}}_{\mat{G}}
\underbrace{
    \begin{pmatrix}
    \theta_1 \\
    \theta_2 \\
    \vdots \\
    \theta_K \\
    \end{pmatrix}}_{\boldsymbol\theta} 
    = \mat{0},
\end{equation}
which is simply linear regression. In practice, we apply singular value decomposition to $\mat{G}=\mat{U}\mat{\Sigma}\mat{V}^T$, where $\mat{U}\in\mathbb{R}^{P\times P}$ and $\mat{V}\in\mathbb{R}^{K\times K}$ are orthogonal matrices, $\mat{\Sigma}\in\mathbb{R}^{P\times K}$ is diagonal with singular values $0\leq\sigma_1\leq\sigma_2\leq\cdots$. We count $\sigma_i$ as effectively zero if $\sigma_i<\epsilon\equiv 10^{-8}$. The number of vanishing singular values, denoted $M$, is equal to the dimensionality of the solution space (null space), which is spanned by the first $M$ columns of $\mat{V}^T$, denoted ${\boldsymbol\Theta}^{(1)}\equiv ({\boldsymbol\theta}_1^{(1)},{\boldsymbol\theta}_2^{(1)},\cdots,{\boldsymbol\theta}_M^{(1)})\in\mathbb{R}^{K\times M}$. The linear structure obviously gives {\bf completeness} (in the space spanned by basis functions), since any solution ${\boldsymbol\theta}$ can be expressed as a linear combination of columns of ${\boldsymbol\Theta}^{(1)}$.

\subsection{Interpretability}
In order to gain more interpretability, we want $\boldsymbol\Theta^{(1)}$ to be sparse. Note that if $\mat{R}\in\mathbb{R}^{M\times M}$ is an orthogonal matrix, the columns of $\boldsymbol\Theta^{(2)}=\boldsymbol\Theta^{(1)} \mat{R}$ also form a set of complete and orthogonal solutions. Therefore we can encourage sparsity by finding and applying the orthogonal matrix that minimizes the following:
\begin{equation}\label{eq:sparsification}
\mat{R}^*=\underset{ \mat{R}^T\mat{R}=\mat{I}}{\rm argmin}\ ||{\boldsymbol\Theta^{(1)}}\mat{R}||_1,\quad 
{\boldsymbol\Theta^{(2)}}= {\boldsymbol\Theta^{(1)}}\mat{R}^*,
\end{equation}
where $||\mat{M}||_1\equiv\sum_{ij}|M_{ij}|$ denotes the $L_1$-norm of a matrix $\mat{M}$, encouraging sparsity.

\subsection{Independence}
Although columns of $\boldsymbol\Theta^{(2)}$ are linearly independent, $H_{\boldsymbol\Theta^{(2)}}$ are not guaranteed to be functionally independent. Take the 1D harmonic oscillator $\mat{x}=(x,p)$, for example. Restricting basis functions to be polynomials in $x$ and $p$ up to the 4th order, there are two solutions:
\begin{equation}
    H_{\boldsymbol\theta_1} = x^2 + p^2, H_{\boldsymbol\theta_2}= H_{\boldsymbol\theta_1}^2 = x^4 + 2x^2p^2 + p^2,
\end{equation}
where $\boldsymbol\theta_1$ and $\boldsymbol\theta_2$ are orthogonal (hence independent), but $H_{\boldsymbol\theta_2}= H_{\boldsymbol\theta_1}^2$, so they are not functionally independent. Consequently, we want a subset of $\boldsymbol\Theta^{(2)}$, denoted $\boldsymbol\Theta^{(3)}$, such that $H_{\boldsymbol\Theta^{(3)}}$ is both independent and complete (i.e., can generate $H_{\boldsymbol\Theta^{(2)}}$). The first question is: how many elements, denoted $c$, does $\boldsymbol\Theta^{(3)}$ have?
As shown in~\cite{poincare2}, $c$ is equal to the rank of the following matrix:
\begin{equation}
    \mat{A}
    =
    \begin{pmatrix}
        \frac{\partial H_{\boldsymbol\theta_1}}{\partial x_1} & \frac{\partial H_{\boldsymbol\theta_2}}{\partial x_1} & \cdots & \frac{\partial H_{\boldsymbol\theta_M}}{\partial x_1} \\
        \frac{\partial H_{\boldsymbol\theta_1}}{\partial x_2} & \frac{\partial H_{\boldsymbol\theta_2}}{\partial x_2} & \cdots & \frac{\partial H_{\boldsymbol\theta_M}}{\partial x_2} \\
        \vdots & \vdots & & \vdots \\
        \frac{\partial H_{\boldsymbol\theta_1}}{\partial x_d} & \frac{\partial H_{\boldsymbol\theta_2}}{\partial x_d} & \cdots & \frac{\partial H_{\boldsymbol\theta_M}}{\partial x_d} \\
    \end{pmatrix}
\end{equation}
which hinges on the fact that gradients of functionally dependent functions are linearly dependent~\footnote{Suppose $H_1$ is dependent on $H_2, H_3$, then there exists a function $g$ such that $H_1=g(H_2,H_3)$. Taking gradients on both es gives $\nabla H_1, \frac{\partial g}{\partial H_2}\nabla H_2+\frac{\partial g}{\partial H_3}\nabla H_3$, which means $\nabla H_1$ is linear combination of $\nabla H_2$ and $\nabla H_3$. So $\mat{A}=[\nabla H_1, \nabla H_2,\nabla H_3]$ has rank two (if $H_2$ and $H_3$ are independent)}. In practice, applying singular value decomposition to $\mat{A}$ gives $\mat{A}=\mat{U}'\mat{\Sigma}'\mat{V}'^T$, where $\mat{U}'\in\mathbb{R}^{d\times d}$ and $\mat{V}'\in\mathbb{R}^{M\times M}$ are orthogonal matrices, and $\mat{\Sigma}'$ is a diagonal matrix with singular values $ s_1\geq s_2\geq \cdots\geq 0$. We count $s_i$ as effectively non-zero if $s_i>\epsilon=10^{-8}$. The number of non-zero singular values is equal to ${\rm rank}(\mat{A})$, which is in turn equal to $c$. %Note that $\mat{A}'\equiv\mat{A}\mat{V}^T=\mat{U}^T\mat{\Sigma}$ is independent in first $c$ columns (all the remaining columns are zero), so first $c$ columns of $\mat{V}^T$ defines a transformation from $\{{\boldsymbol\theta_i}, 1\leq i\leq M\}$ to a set of new parameters $\{{\boldsymbol\theta_j'}, 1\leq j\leq c\}$ such that $\{H_{\boldsymbol\theta_j'}(\mat{x}), 1\leq j\leq c\}$ are functionally independent.
After determining $c$, we aim to obtain $\boldsymbol\Theta^{(3)}$ by selecting $c$ elements from $\boldsymbol\Theta^{(2)}$. The selection process is as follows: (1) We assign each conserved quantity a complexity score (based on entropy~\footnote{For a vector ${\boldsymbol\theta}$, we can associate it to a  probability distribution $p_i\equiv|\theta_i|/(\sum_i|\theta_j|)$, whose entropy is $S=-\sum_i p_i{\rm log}p_i$.}) and sort them from the simplest to the most complex. (2) Starting from an empty set $\boldsymbol\Theta^{(3)}$, looping over element $\boldsymbol\theta \in\boldsymbol\Theta^{(2)}$, we add $\boldsymbol\theta$ to $\boldsymbol\Theta^{(3)}$ if $H_{\boldsymbol\theta}$ is independent of $H_{\boldsymbol\Theta^{(3)}}$ (functions already added), until $\boldsymbol\Theta^{(3)}$ contains $c$ elements.

\section{Results}

To better illustrate SID, we apply it to three dynamical systems.

\begin{figure}[tbp]

\begin{subfigure}{0.2\textwidth}
\includegraphics[width=1.0\linewidth]{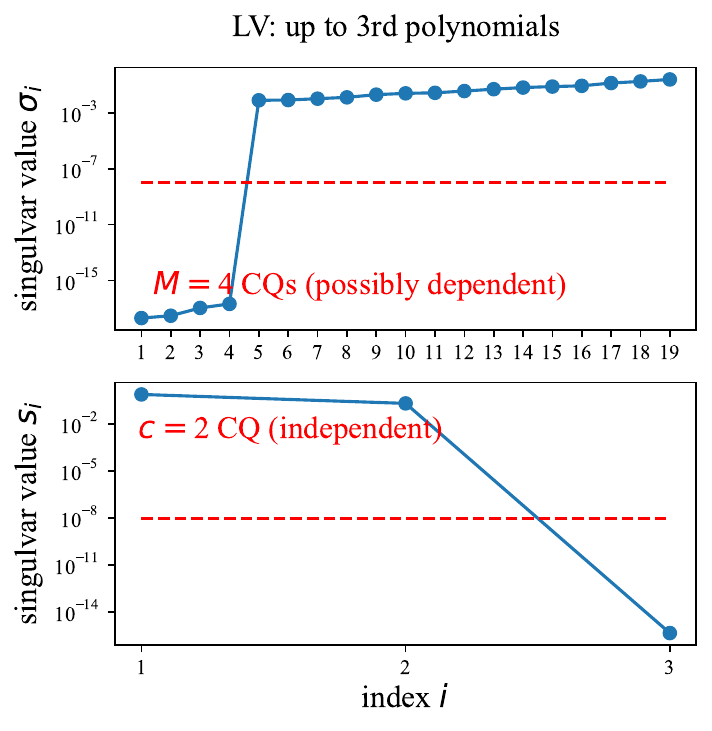} 
\caption{}
\label{fig:LV_sv}
\end{subfigure}
\begin{subfigure}{0.27\textwidth}
\includegraphics[width=1.0\linewidth]{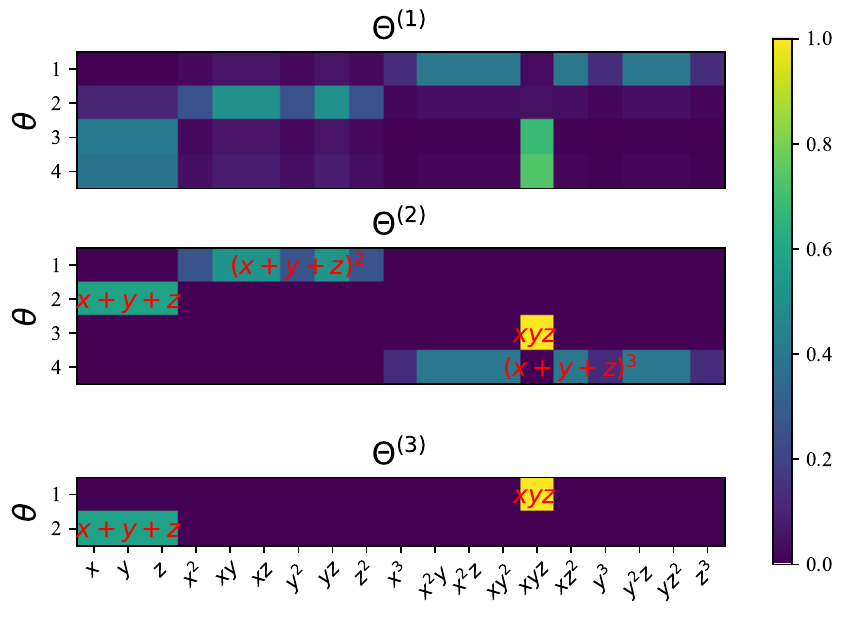}
\caption{}
\label{fig:LV_Theta}
\end{subfigure}
\caption{Three-species Lokta-Volterra equation. SID correctly discovers that: (a) there are $4$ CQs (top) in polynomials up to 3rd order, yet only 2 of them (bottom) are independent. (b) coefficients of conserved quantities $\boldsymbol\Theta^{(i)} (i=1,2,3)$, with more interpretability.}
\label{fig:LV}
\end{figure}

\begin{figure*}[htbp]
    \centering
    \includegraphics[width=0.9\linewidth]{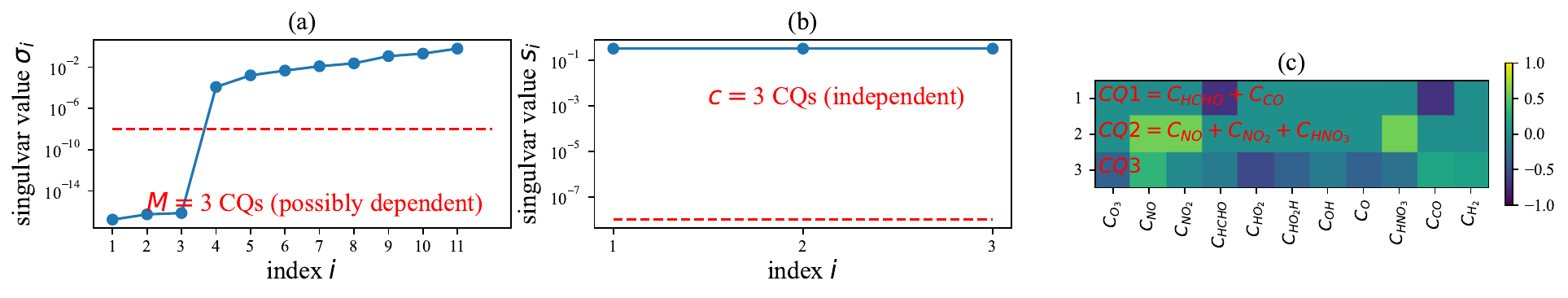}
    \caption{(a)(b) SID discovers 3 independent conserved quantities in the ozone photochemical production model. (c) The coefficients of these linear conserved quantities. The first two correspond to the known carbon and nitrogen conservation, while $CQ_3$ is identified for the first time.}
    \label{fig:chem}
\end{figure*}

\begin{figure}[htbp]
    \centering
    \includegraphics[width=0.8\linewidth]{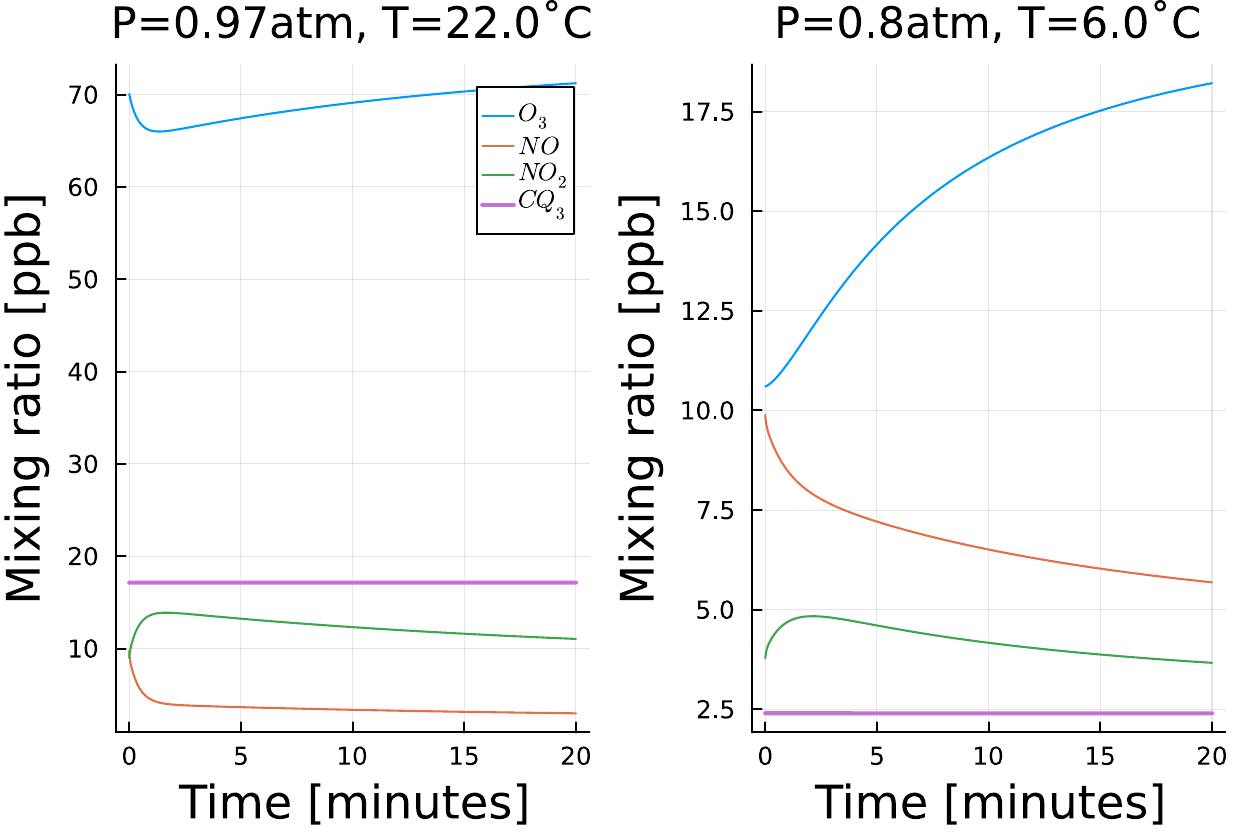}
    \caption{The evolution of  concentrations in simulation. Under both conditions, $CQ_3$ is well conserved.}
    \label{fig:CQ34}
\end{figure}

\subsection{Systems biology}  Our first application is from  systems biology. The Lotka–Volterra equations (LV hereafter) describe how population of  many species evolve in time via inter-species interactions. We study this particular equation: 
\begin{equation}
    \dot{x}=x(y-z), \dot{y}=y(z-x),\dot{z}=z(x-y)
\end{equation}
with two known conserved quantities $H_1 = x+y+z$ and $H_2=xyz$. We define basis functions to be all polynomials up to 3rd order (including $K=19$ terms, shown in FIG.~\ref{fig:LV_Theta}). We draw $P=100$ data points from the standard Gaussian distribution, i.e., $\mat{x}_i\equiv (x_i,y_i,z_i)\sim \mathcal{N}(0,\mat{I}_{3\times 3})\ (1\leq i\leq P)$. 

Within the function space spanned by the basis functions, $M=4$ conserved quantities are found, since FIG.~\ref{fig:LV_sv} top shows that there are 4 vanishing singular values $\sigma_i$. The coefficients of four CQs are somewhat mixed, shown in FIG.~\ref{fig:LV_Theta} top. After sparsification (Eq.~(\ref{eq:sparsification})), the coefficients become less entangled, shown in FIG.~\ref{fig:LV_Theta} middle, although their represented conserved quantities are still dependent. Among the 4 CQs, only $c=2$ are independent, since FIG.~\ref{fig:LV_sv} bottom shows that there are 2 non-vanishing singular values $s_i$. FIG.~\ref{fig:LV_Theta} bottom shows the final outputs: the two conserved quantities agree with our prior knowledge.

What will happen if we choose the set of basis functions to be smaller or larger? (1) smaller: if we include polynomials up to the 1st or the 2nd order ($K=3$ or $K=9$) only, then only $H_1$ can be discovered by SID, while $H_2$ is not discovered. (2) larger: if we instead include polynomials up to order 4,5,6 ($K=34, 55, 83$), then both $H_1$ and $H_2$ are still discovered, although the sparsification and independence process may take longer than with only 3rd order polynomials. Detailed results are included in Appendix~\ref{app:LV_results}. The take-home message is that SID does not replace human scientists since it requires the input of basis functions (formulate hypothesis) from human scientists. %However, SID is a helpful assistant such that when human scientists interact with it, the discovery process is facilitated. Put it another way, 
SID is good at testing hypotheses (which could be technical and tedious), however, it is human scientists who formulate hypotheses (which requires creativity).

% \zm{We have tried our best to interpret $CQ_3$ in terms of known facts, but these attempts have all failed. In particular, we suspected $CQ_3$ is due to approximate conservation of Hydrogen atoms, but it turned out not to be case (see Appendix \ref{sec:CQ3_not_H}). This implies that $CQ_3$ might be a highly non-trivial conserved quantity, which is worth thorough study in future works.}

% Typical simulation trajectories are shown in Figure~\ref{fig:CQ34}, where $CQ_3$ is not only conserved under the training condition (left), but also holds in a different condition (right).  
% A caveat is that these quantities are not always conserved in the 1000 validation simulations: the maximum coefficient of variation is 50\% for $CQ3$ and 113\% for $CQ4$. This means that  

\subsection{Fluid mechanics}

Arguably the biggest puzzle in fluid mechanics is turbulence~\cite{falkovich2006lessons, Constantin2007OnTE}. Turbulence, and chaos in general, are due to lack of sufficient conserved quantities. Therefore, studying conserved quantities of fluid systems is relevant to understanding turbulence. As a preliminary step, we study conserved quantities of a fluid element in ideal fluid (zero viscosity and incompressible). In 2D (3D), The fluid element is a triangle (tetrahedron), which is represented by its 3 (4) vertices~\cite{pumir2013tetrahedron}. Effectively, we can view the system as 3 (4) ``free" particles, with the only constraint being that the area (volume) of the triangle (tetrahedron) should remain unchanged. The equations of motion are included in Appendix \ref{app:fluid}, which appear a bit intimidating (especially for 3D).

Fluid dynamics experts (including some authors of this paper) have attempted to find the conserved quantities with pencil and paper. Given the complexity of the calculations, it is impressive that they found 8 (12) conserved quantities for 2D (3D). However, they were unsure whether there were more undiscovered conserved quantities, and whether the discovered ones were in their simplest. So we turn to SID for help. The results below are for the basis function set selected to be polynomials up to 2nd (3rd) order for 2D (3D), but more polynomial orders are also tried in Appendix~\ref{app:fluid}.

For the 2D case, SID finds 8 conserved quantities, agreeing with experts' expectation. Interestingly, the conserved quantities found by SID appear to be simpler. In fact, all the conserved quantities discovered by SID are 1st or 2nd order polynomials, while experts found a 4th order polynomial, which we find to b a combination of two 2nd order are conserved quantities discovered by SID. 

For the 3D case, SID finds 14 conserved quantities, while experts found only 12. The two new conserved quantities can be interpreted as the angular momentum in the center of mass (COM) frame. They are non-trivial because it is easy to (falsely) think that the COM angular momentum is dependent on the angular momentum and the linear momentum~\footnote{Also, the COM angular momentum has three components in 3D, so why are there only two new conserved quantities, not three? It turns out that the third component is not independent of the 14 conserved quantities.}. %(hence not independent), but the intuition 
%is wrong since it neglects cross terms. (2) The COM angular momentum has three components in 3D, so why are there only two more conserved quantities not three? It turns out that the third component is not independent of the 14 conserved quantities. 
Although humans alone will probably get the results right at the end of the day without SID, SID can take care of subtle details automatically, thus saving human experts' mental labor to a great extent.

\subsection{Atmospheric chemistry}

We next apply SID to a truncated atmospheric chemistry model of photochemical ozone production \cite{StuWex2022conservation}, where an exotic new conserved quantity is found. This simplified dynamical system contains 11 species and 10 reactions involved in ozone formation, including NOx, organic, and radical chemistry \cite{StuWex2020mass}. A key characteristic of this system is conservation of carbon and nitrogen atoms, $H_C$ and $H_N$, respectively. Though species in this model contain two other elements, hydrogen and oxygen, neither are conserved, as $H_2O$ molecules are not one of the 11 species whose concentrations are tracked, and diatomic oxygen $O_2$ is treated as an infinite source and sink due to its abundance. $H_C$ and $H_N$ are implied in the coefficients of a stoichiometric matrix $\mat{B} \in \mathbb{Z}^{11,10}$ used in prior work to enforce conservation of atoms in machine learning surrogate models \cite{StuWex2022conservation}. $H_C$ and $H_N$ can be represented by linear combinations of species concentrations, the coefficients of which form a basis for the null space of $\mat{B}^T$. Further details are provided in Appendix \ref{app:chem}.

We applied SID to simulation trajectories expecting to discover up to 2 conserved quantities which are linear combinations of concentrations. The training data are points on simulation trajectories at pressure $P=0.95$ atm and temperature $T=20.0^\circ $C. 
%SID surprisingly discovers 3 conserved quantities, as shown in Figure~\ref{fig:chem}. 
As shown in FIG.~\ref{fig:chem}, besides $H_C$ and $H_N$, SID surprisingly discovers a third conserved quantity $CQ_3$ that is a linear combination of species concentrations ($C_X$ means the concentration of $X$):
\begin{equation}
\begin{aligned}
    CQ_3 \approx  &\ 6C_{O_3} - 5 C_{NO}+C_{NO_2} + 3C_{HCHO} \\
    & + 9C_{HO_2} + 6C_{HO_2H} + 2C_{OH} + 6C_{O}  \\
    & + 4C_{HNO_3} - 3C_{CO} - 2.21C_{H_2}
\end{aligned}
\end{equation}
This additional quantity is linearly independent of $H_C$ and $H_N$ and is not in the null space of $\mat{B}^T$. $CQ_3$ has a relative variation of less than 0.1\% in 995 of 1000 simulated cases. Two representative simulation trajectories are shown in Figure~\ref{fig:CQ34}, where $CQ_3$ holds under different chemical and meteorological conditions.  The evolving concentrations of $O_3$, $NO$, and $NO_2$ are included as contrasts to the invariance of $CQ_3$. We have not yet identified the underlying cause of $CQ_3$, and whether it is physically exact or numerically approximate.
We have ruled out symmetry corresponding to hydrogen conservation: when explicitly incorporating production of $H_2O$ as an additional buildup species, SID identifies approximate hydrogen conservation as well as a fourth conserved quantity  (see Appendix \ref{sec:CQ3_not_H}). This implies that $CQ_3$ might be a non-trivial conserved quantity that is worth thorough study in future work.

\section{Conclusions}
We have presented an algorithm SID to automatically discover conserved quantities from dynamical equations. In constrast to previous blackbox models, SID is guaranteed to be robust and interpretable thanks to its algorithmic simplicity. We demonstrate the power of SID on two examples in atmospheric chemistry and fluid mechanics, revealing new  conserved quantities hitherto unknown to human experts. Although SID does not replace human scientists, it is a helpful assistant that can facilitate the discovery process. Promising future directions include applying SID to a broader range of applications and explicitly deal with symmetries that users may want to impose. %So far we do not have any constraint on the symmetry of conserved quantities. It might also be useful to provide an interface to users where they can specify symmetry constraints. 

\bibliography{poincare_application}

\newpage

\onecolumngrid

\appendix

\section{More Lotka-Volterra Results}\label{app:LV_results}
\begin{figure}[htbp]
\includegraphics[width=0.16\linewidth]{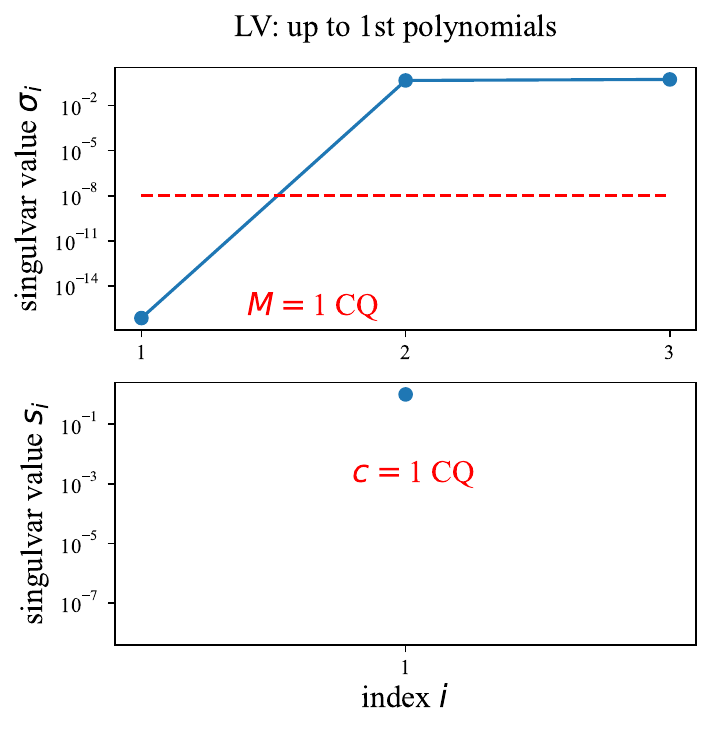} 
\includegraphics[width=0.16\linewidth]{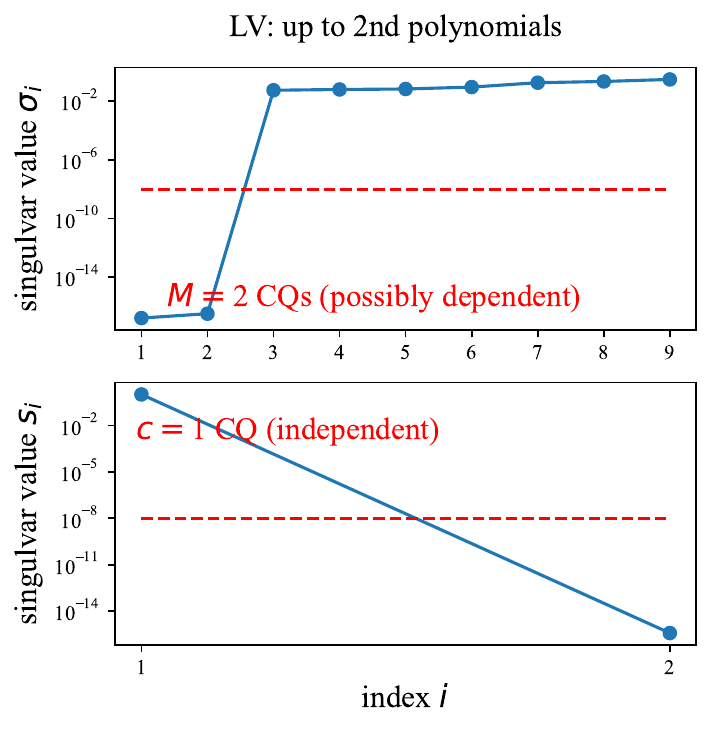} 
\includegraphics[width=0.16\linewidth]{./fig/LV_sv_poly3.pdf} 
\includegraphics[width=0.16\linewidth]{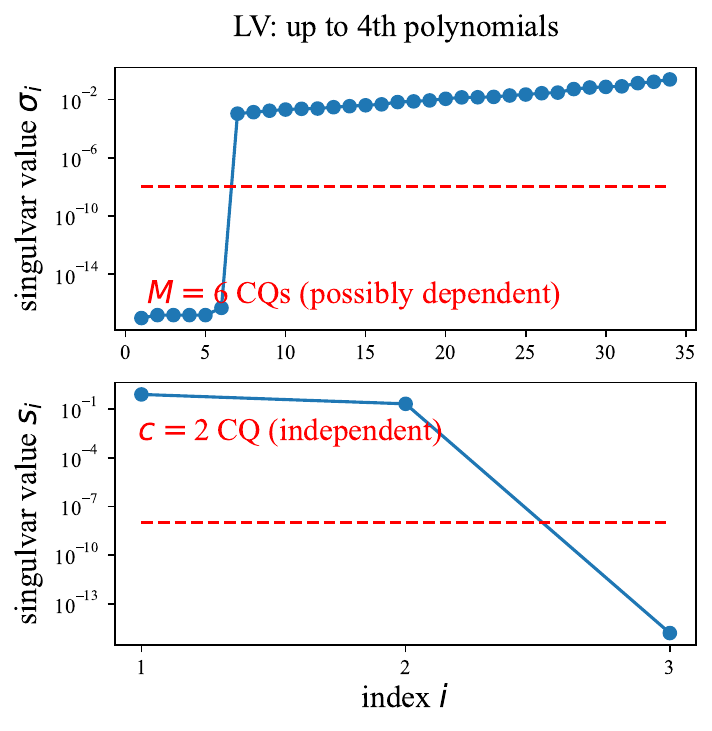} 
\includegraphics[width=0.16\linewidth]{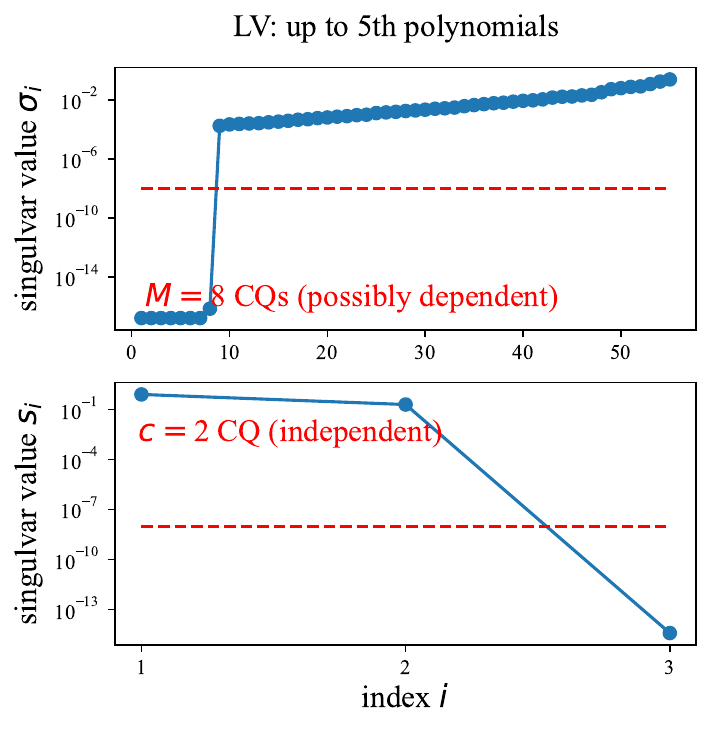} 
\includegraphics[width=0.16\linewidth]{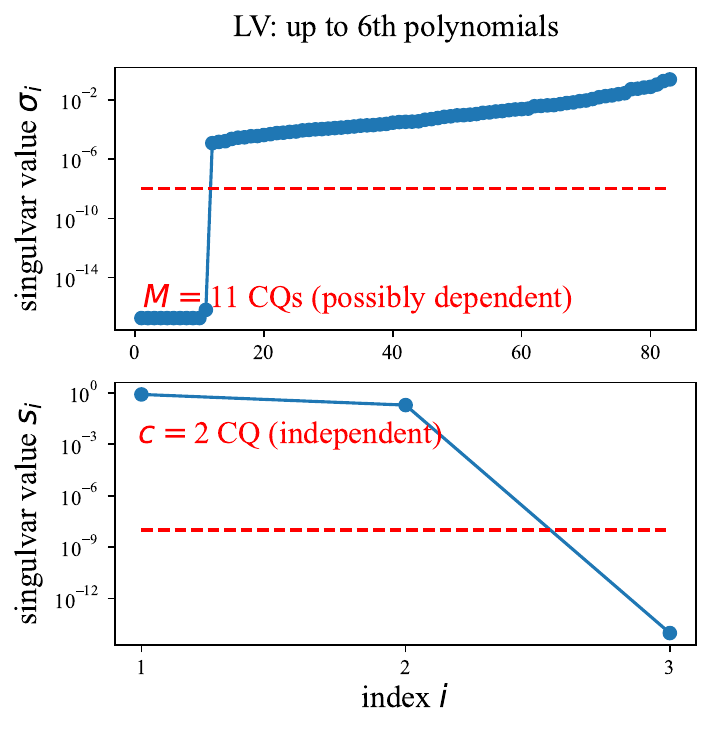} 
\caption{SID's results of the three-species Lokta-Volterra equation with different selections of basis functions. The number of conserved quantities is correct (2) if the basis set is large enough (polynomial order $n\geq 3$).}
\label{fig:LV_sweep_order}
\end{figure}

We study the dependence of SID's outputs on the selection of basis functions. For the Lotka-Volterra example, there are two conserved quantities $H_1=x+y+z$ (1st order polynomial) and $H_2=xyz$ (3rd order polynomial). We choose the set of basis functions to be polynomials up to $n=\{1,2,3,4,5,6\}$ order. As shown in Figure~\ref{fig:LV_sweep_order}, we find that when $n=1$ or $n=2$, only $H_1$ is discovered. For $n=3,4,5,6$, both $H_1$ and $H_2$ are discovered by SID.

\section{Atmospheric chemistry}\label{app:chem}

\subsection{Background}

The simplified model of ozone photochemistry is written in the Julia computing language and represents the broad temporal scales and variability in the prediction of ozone and its precursors. This model relates 11 chemical species via 10 reactions \cite{StuWex2022conservation}. The state of the system is described by the concentration vector
\begin{equation}
    \mat{x} \equiv [C_{O_3}, C_{NO}, C_{NO_2}, C_{HCHO}, C_{HO_2}, C_{HO_2H}, C_{OH}, C_{O}, C_{HNO_3}, C_{CO}, C_{H_2}]\equiv [C_1,C_2,\dots].
\end{equation}
Species and reactions form two distinct sets in a bipartite network \cite{silva2021graph} and can be related via a weighted biadjacency matrix $\mat{B} \in \mathbb{Z}^{11,10}$ that specifies the stoichiometry of the chemical system. Conservation of carbon and nitrogen atoms are implicit in the coefficients of $\mat{B}$.

\setcounter{MaxMatrixCols}{11}
\begin{equation}\label{eq:stoichiom-matrix}
\mat{B}=\begin{bmatrix}[c|rrrrrrrrrr]
&R1&R2&R3&R4&R5&R6&R7&R8&R9&R10 \\\hline O_3&0&1&-1&0&0&0&0&0&0&0\\NO&1&0&-1&0&0&0&-1&0&0&0 \\{NO}_2&-1&0&1&0&0&0&1&-1&0&0\\HCHO&0&0&0&-1&-1&-1&0&0&0&0\\{HO}_2&0&0&0&2&0&1&-1&0&0&1\\{HO}_2H&0&0&0&0&0&0&0&0&-1&-1\\OH&0&0&0&0&0&-1&1&-1&2&-1\\O&1&-1&0&0&0&0&0&0&0&0\\HNO_3&0&0&0&0&0&0&0&1&0&0\\CO&0&0&0&1&1&1&0&0&0&0\\H_2&0&0&0&0&1&0&0&0&0&0\\\end{bmatrix}
\end{equation}

A system of coupled ODEs represents the time evolution of the concentrations of species.  In this system, each reaction $j$ has a rate $r_j$ defined by the law of mass action %[note from Obin: we can remove this reference if not appropriate/necessary] % Ziming: This looks fine.

\begin{equation}\label{eq:mass-action}
r_j = k_j \prod_{i|B_{ij}<0} C_i^{\abs{B_{ij}}}
\end{equation}
where the reaction rate constant $k_j$ is determined by an Arrhenius relation. The rate of change of concentration for each species $i$ is the sum of reaction rates describing its production and loss:

\begin{equation}\label{eq:concentration-rate}
\frac{dC_i}{dt} =\sum_j B_{ij}r_j,
\end{equation}

Substituting Eq.~(\ref{eq:mass-action}) for $r_j$ in Eq.~(\ref{eq:concentration-rate}), the full set of ODEs for the system in Eq.~(\ref{eq:stoichiom-matrix}) can be written as

\begin{equation}\label{eq:chemistry-ode-system}
    \begin{aligned}
        \frac{dC_{O_3}}{dt}&=k_2C_O\ -{\ k}_3C_{O_3}C_{NO}, \\
        \frac{dC_{NO}}{dt}&=k_1C_{NO_2}\ -{\ k}_3C_{O_3}C_{NO}-{\ k}_7{C_{NO}C}_{HO_2},\\
        \frac{dC_{NO_2}}{dt}&={-\ k}_1C_{NO_2}\ +{\ k}_3C_{O_3}C_{NO}+k_7C_{NO}C_{HO_2}-{\ k}_8C_{NO_2}C_{OH}, \\
        \frac{dC_{HCHO}}{dt}&={-\ k}_4C_{HCHO}\ -{\ k}_5C_{HCHO}-{\ k}_6C_{HCHO}C_{OH}, \\
        \frac{dC_{HO_2}}{dt}&={2k}_4C_{HCHO}+\ {\ k}_6C_{HCHO}C_{OH}-{\ k}_7{C_{NO}C}_{HO_2}+k_{10}C_{HO_2H}C_{OH},\\
        \frac{dC_{HO_2H}}{dt}&=-k_9C_{HO_2H}-k_{10}C_{HO_2H}C_{OH}, \\
        \frac{dC_{OH}}{dt}&=-{\ k}_6C_{HCHO}C_{OH}+{\ k}_7{C_{NO}C}_{HO2}-{\ k}_8C_{NO2}C_{OH}+{2k}_9C_{HO2H}-k_{10}C_{HO2H}C_{OH}, \\
        \frac{dC_O}{dt}&=k_1C_{NO2}-k_2C_O, \\
        \frac{dC_{HNO3}}{dt}&={\ k}_8C_{NO2}C_{OH}, \\
        \frac{dC_{CO}}{dt}&={\ k}_4C_{HCHO}\ +{\ k}_5C_{HCHO}+{\ k}_6C_{HCHO}C_{OH}, \\
        \frac{dC_{H2}}{dt}&={\ k}_5C_{HCHO}.
    \end{aligned}
\end{equation}

The table below gives example values for $k$ at T=298K and P=1 atm, with photolytic constants representative of mid-day conditions in Southern California. In the model and the testing of the conserved quantities, rate constants are dependent on reaction-specific parameters, as well as temperature and pressure, determined by a dry adiabatic lapse rate and the hypsometric equation whose values are determined by an altitude $z$ drawn from a uniform distribution on the interval $[0,2000]$ meters. 

\begin{table}[H]
\centering
\begin{tabular}{|c|c|c|c|c|c|c|c|c|c|c|}
\cline{1-11}
  rate constant & $k_1$  & $k_2$ & $k_3$ & $k_4$ & $k_5$ & $k_6$ & $k_7$ & $k_8$ & $k_9$ & $k_{10}$ \\ \cline{1-11}
  value & 0.5 & 22.179 & 26.937 & 0.015 & 0.022 & 13844.97 & 12652.43 & 15454.98 & 0.0003 & 2492.71  \\ \cline{1-11}
\end{tabular}
\end{table}

To reduce the numerical stiffness of this system, the pseudo-steady state approximation is applied to two species (hydroxyl radical $OH$ and atomic oxygen $O$) whose rates of change are approximated as zero and concentrations rewritten as algebraic expressions of the other species concentrations.  The concentration vector $C$ is restricted to its semi-positive half space (negative concentrations are unphysical) initialized from uniform distributions where ozone is selected to range from 0-100 ppb and other ranges informed by initialization ranges in prior work \cite{Kelp2020}.

\subsection{Conserved quantity results}

{\bf Known conserved quantities} We expect $H_C=C_{HCHO}+C_{CO}$ and $H_N=C_{NO}+C_{NO_2}+C_{HNO_3}$ to be conserved, as this is built into the model. %These can be calculated using linear combinations of species, specifically, summing all of the carbon and nitrogen atoms. Oxygen and hydrogen atoms are not conserved, as diatomic oxygen $O_2$ is modeled as an infinite source and sink due to its abundance in the atmosphere and concentration of water $H_2O$ is not tracked. Though a pseudo steady state approximation is made for two other species, $OH$ and $O$, where their instantaenous rates of change are approximated as zero, their concentrations are not constant across timesteps, but rather specified by nonlinear, algebraic expressions of other concentrations in the model.

{\bf SID} SID identifies three conserved quantities, two of which are $H_C$ and $H_N$. To validate the results of SID, we initialize one thousand 20-minute simulations under the above conditions. The known conservation properties hold in the model to very high precision, with coefficients of variation (mean-normalized standard deviation) of less than $10^{-14}$ for all simulations.  The rediscovered conservation laws for carbon and nitrogen atom closure also display high accuracy, with maximum coefficients of variation over all 1000 cases on the order of $10^{-11}$ and  $10^{-13}$ respectively. $CQ_3$ is not conserved to such high precision in the 1000 validation simulations: its maximum coefficient of variation is $0.96\%$.  However, the coefficient of variation of $CQ_3$ is within 0.1\% in all but 5 of 1000 cases, and has a 95 percentile over all cases of 0.02\%.

\setcounter{MaxMatrixCols}{20}
The coefficients discovered by SID are:
\begin{equation}
    \mat{B} = 
    \begin{pmatrix}[c|ccccccccccc]
     & O_3 & NO & NO_2 & HCHO & HO_2 & HO_2H & OH & O & HNO_3 & CO & H_2 \\\hline
    CQ_1 & 0 & 0 & 0 &  0.707 & 0 & 0 & 0 & 0 & 0 & 0.707 & 0  \\
    CQ_2 & 0 & 0.577 & 0.577 & 0 & 0 & 0 & 0 & 0 & 0.577 & 0 & 0 \\
    CQ_3 & 0.370 & -0.310 & 0.061 & 0.185 & 0.555 & 0.370 & 0.185 & 0.370 & 0.247 & -0.185 & -0.135 \\
    \end{pmatrix}
\end{equation}

Note that the ratios of non-zero coefficients are very close to rational numbers, so if we divide rows by $[0.707,0.577,0.062]$, we approximately get
\begin{equation}
    \mat{B}' = 
    \begin{pmatrix}[c|ccccccccccc]
     & O_3 & NO & NO_2 & HCHO & HO_2 & HO_2H & OH & O & HNO_3 & CO & H_2 \\\hline
    CQ_1 & 0 & 0 & 0 &  1 & 0 & 0 & 0 & 0 & 0 & 1 & 0  \\
    CQ_2 & 0 & 1 & 1 & 0 & 0 & 0 & 0 & 0 & 1 & 0 & 0 \\
    CQ_3 & 6 & -5 & 1 & 3 & 9 & 6 & 3 & 6 & 4 & -3 & -2 \\
    \end{pmatrix}
\end{equation}
The snapped $CQ_3$ is not quite as conserved: in 910 out of 1000 cases, $CQ_3$ is conserved within 0.1\%. However, if we do not snap the coefficient for $H_2$ as 2 (instead set 2.213115), $CQ_3$ holds to within 0.1\% for 995 out of 1000 cases. 

%We have tried our best to interpret $CQ_3$ in terms of known facts, but these attempts have all failed. In particular, we suspected $CQ_3$ is due to approximate conservation of Hydrogen atoms, but it turned out not to be case (see below).

\subsection{$CQ_3$ is independent of hydrogen atom conservation}\label{sec:CQ3_not_H}

In the 11-species model, the number of hydrogen atoms is not conserved because the net production of $H_2O$ is ignored. The atomic composition of each species is not explicitly represented in the model but rather implied by the stoichiometric coefficients of the reactions, which ensure carbon and nitrogen conservation. We find that $CQ_3$ is balanced in 8 out of 10 individual reactions. This motivates construction of an experiment to assess whether SID could be inferring some incorrect composition of hydrogen in each species (which is conserved as represented by $CQ_3$). To assess this, we augment the model to track net production of $H_2O$ as a 12th species so that the number of hydrogen atoms becomes:
\begin{equation}\label{HH}
    H_H = 2C_{HCHO} + C_{HO_2} + 2C_{HO_2H} + C_{OH} + C_{HNO_3} + 2C_{H_2} + 2C_{H_2O}.
\end{equation}
%However, since SID does not distinguish between exact or approximate conserved quantities (as long as the conservation loss is below $\epsilon$), $CQ_3$ could still be related to the number of hydrogen atoms. The number of hydrogen atoms for the 11-species model is
%\begin{equation}
%    H_H = 2C_{HCHO} + C_{HO_2} + 2C_{HO_2H} + C_{OH} + C_{HNO_3} + 2C_{H_2}.
%\end{equation}
%However, it is clear that $H_H$ is not a linear combination of $CQ_1$, $CQ_2$ and $CQ_3$. 

We note that (\ref{HH}) is not fully conserved in the 12-species model due to the pseudo-steady state approximation of $OH$, but is approximately conserved, with a maximum mean-normalized standard deviation of 4.2\% over all cases, and holding in 972 out of 1000 cases to within 0.1\%. With $H_2O$ represented in the 12-species model, there are three known conserved quantities: conservation of carbon and nitrogen and approximate conservation of hydrogen.  If $CQ_3$ indeed arises from an incorrectly inferred composition of hydrogen in the set of 11 species, we would expect three conserved quantities for the 12-species model, with the third quantity representing the known hydrogen balance in (\ref{HH}). However, SID finds 4 conserved quantities, demonstrating that there is indeed an additional conserved quantity independent of hydrogen atom conservation.

\section{Fluid mechanics}\label{app:fluid}

\subsection{Background}

\begin{figure}[htbp]
    \centering
    \includegraphics[width=0.8\linewidth, trim=0cm 0cm 0cm 0cm]{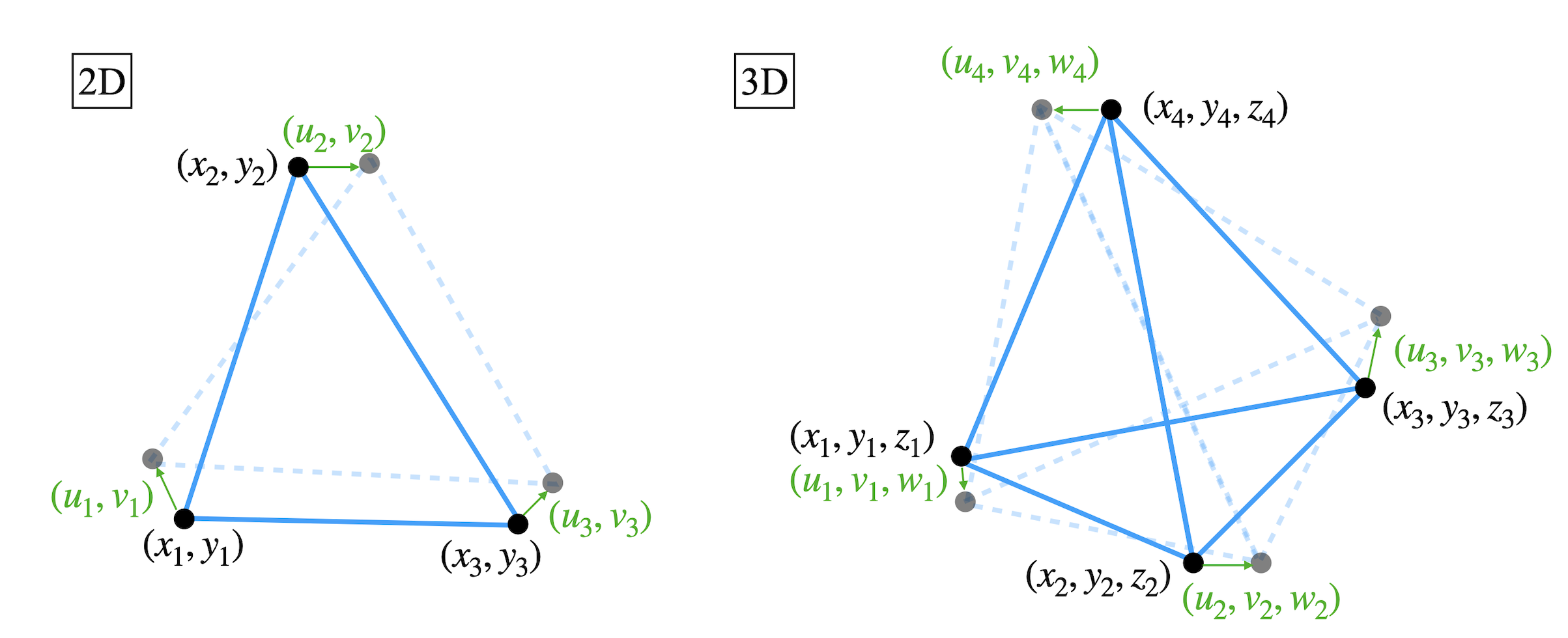}
    \caption{Illustration of fluid element motion. Left: 2D, right: 3D.}
    \label{fig:fluid_illustration}
\end{figure}

Fluid mechanics systems can be very complicated, especially when exhibiting turbulence~\cite{falkovich2006lessons,Constantin2007OnTE}. Turbulence, and chaos in general, emerge due to lack of sufficient conserved quantities. To understand turbulence, it is thus relevant to study conserved quantities of fluid systems. We will deal with ideal fluid, and partition the continuum of fluid into small elements (triangles in 2D or irregular tetrahedra in 3D)~\cite{pumir2013tetrahedron}. For ideal incompressible fluid, the only internal interaction between elements is via pressure, which may perform work on individual element (pressure gradients as external forces), but the total work for the whole fluid is zero. Since we care about global conserved quantities of the whole fluid, without loss of generality, we can assume zero external forces on each individual element. Therefore, we can first study local conserved quantities (for each element), and the global conserved quantities of the whole fluid is the summation (integral) of conserved quantities of individual elements. This additivity property allows us to study a fluid element first. A fluid element is described by its vertices (3 in 2D, 4 in 3D), illustrated in Figure~\ref{fig:fluid_illustration}. Every vertex $i$ is in turn described by its position $(x_i,y_i,z_i)$ and velocity $(u_i,v_i,w_i)\equiv (\dot{x}_i,\dot{y}_i,\dot{z}_i)$. The motion of the fluid element can thus be described by the motion of its vertices. Because the fluid element is small, one can imagine unit masses put on the vertices, with no mass elsewhere. So vertices move like free particles, with the only interaction (constraint) being that the area/volume is conserved. Note that since we partition the fluid to tetrahedra (which lack knots or linkages), we may miss topological invariants (e.g., helicity, which is the dot product of the velocity and vorticity).

\subsubsection{2D} Next, we use Lagrangian mechanics to obtain equations of motion in 2D, which is straightforwardly generalized to 3D. The phase space is 12D: $\mat{x}=(x_1,y_1,u_1,v_1,x_2,y_2,u_2,v_2,x_3,y_3,u_3,v_3)$. The area of the triangle is
\begin{equation}
    A = x_1y_2 + x_2y_3 + x_3y_1 - x_2y_1 - x_3y_2 - x_1y_3.
\end{equation}
Since the area is conserved (we are assuming an incompressible fluid), we impose that 
\begin{equation}\label{eq:area_constraint}
    \frac{dA}{dt} = u_1(y_2 - y_3) + v_1(x_3 - x_2) + u_2(y_3 - y_1) + v_2(x_1 - x_3) + u_3(y_1 - y_2) + v_3(x_2 - x_1) = 0
\end{equation}
The Lagrangian is the kinetic energy of three particles, plus the area constraint ($\lambda$ is the Lagrange multiplier):
\begin{equation}
    \mathcal{L} = \frac{1}{2}(u_1^2+v_1^2+u_2^2+v_2^2+u_3^2+v_3^2) + \lambda (x_1y_2 + x_2y_3 + x_3y_1 - x_2y_1 - x_3y_2 - x_1y_3)
\end{equation}
The Euler-Lagrange equations give
\begin{equation}
    \begin{aligned}
    &\dot{x_1}=u_1,\quad \dot{u_1}=\lambda(y_2-y_3),\quad
    \dot{y_1}=v_1,\quad \dot{v_1}=\lambda(x_3-x_2), \\
    &\dot{x_2}=u_2,\quad \dot{u_2}=\lambda(y_3-y_1),\quad
    \dot{y_2}=v_2,\quad \dot{v_2}=\lambda(x_1-x_3), \\
    &\dot{x_3}=u_3,\quad \dot{u_3}=\lambda(y_1-y_2),\quad \dot{y_3}=v_3,\quad \dot{v_3}=\lambda(x_2-x_1),
    \end{aligned}
\end{equation}
where $\lambda$ can be solved by noticing $\frac{d^2A}{dt^2}=0$:
\begin{equation}
    \lambda =\frac{(u_1 v_2 - u_1 v_3 - u_2 v_1 + u_2 v_3 + u_3 v_1 - u_3 v_2)}{(-x_1^2 + x_1 x_2 + x_1 x_3 - x_2^2 + x_2 x_3 - x_3^2 - y_1^2 + y_1 y_2 + y_1 y_3 - y_2^2 + y_2 y_3 - y_3^2)}.
\end{equation}

\subsubsection{3D} The derivation of equations of motion is similar to 2D, only with area switched to volume. The phase space is 24D: $\mat{x}=(x_i, y_i, z_i, u_i, v_i, w_i), i=1,2,3,4$. We list the equations here:

\begingroup
\allowdisplaybreaks
\begin{equation}
    \begin{aligned}
    &\dot{x_1}=u_1, \dot{u_1}=\lambda(-y_2 z_3 + y_2 z_4 + y_3 z_2 - y_3 z_4 - y_4 z_2 + y_4 z_3), \\
    &\dot{y_1}=v_1, \dot{v_1}=\lambda(x_2 z_3 - x_2 z_4 - x_3 z_2 + x_3 z_4 + x_4 z_2 - x_4 z_3), \\
    &\dot{z_1}=w_1, \dot{w_1}=\lambda(-x_2 y_3 + x_2 y_4 + x_3 y_2 - x_3 y_4 - x_4 y_2 + x_4 y_3), \\
    &\dot{x_2}=u_2, \dot{u_2}=\lambda(y_1 z_3 - y_1 z_4 - y_3 z_1 + y_3 z_4 + y_4 z_1 - y_4 z_3), \\
    &\dot{y_2}=v_2, \dot{v_2}=\lambda(-x_1 z_3 + x_1 z_4 + x_3 z_1 - x_3 z_4 - x_4 z_1 + x_4 z_3), \\
    &\dot{z_2}=w_2, \dot{w_2}=\lambda(x_1 y_3 - x_1 y_4 - x_3 y_1 + x_3 y_4 + x_4 y_1 - x_4 y_3), \\
    &\dot{x_3}=u_3, \dot{u_3}=\lambda(-y_1 z_2 + y_1 z_4 + y_2 z_1 - y_2 z_4 - y_4 z_1 + y_4 z_2), \\
    &\dot{y_3}=v_3, 
    \dot{v_3}=\lambda(x_1 z_2 - x_1 z_4 - x_2 z_1 + x_2 z_4 + x_4 z_1 - x_4 z_2), \\
    &\dot{z_3}=w_3, 
    \dot{w_3}=\lambda(-x_1 y_2 + x_1 y_4 + x_2 y_1 - x_2 y_4 - x_4 y_1 + x_4 y_2), \\
    &\dot{x_4}=u_4, 
    \dot{u_4}=\lambda\cdot(y_1 z_2 - y_1 z_3 - y_2 z_1 + y_2 z_3 + y_3 z_1 - y_3 z_2), \\
    &\dot{y_4}=v_4,
    \dot{v_4}=\lambda(-x_1 z_2 + x_1 z_3 + x_2 z_1 - x_2 z_3 - x_3 z_1 + x_3 z_2), \\
    &\dot{z_4}=w_4, 
    \dot{w_4}=\lambda(x_1 y_2 - x_1 y_3 - x_2 y_1 + x_2 y_3 + x_3 y_1 - x_3 y_2), \\
    \end{aligned}
\end{equation}
\endgroup
where
\begin{equation}
\begin{aligned}
    \lambda = &\ -2p/q, \\
    p = &\ x_1 (-v_2 w_3 + v_2 w_4 + v_3 w_2 - v_3 w_4 - v_4 w_2 + v_4 w_3) + y_1 (u_2 w_3 - u_2 w_4 - u_3 w_2 + u_3 w_4 + u_4 w_2 - u_4 w_3) \\
    &+ z_1 (-u_2 v_3 + u_2 v_4 + u_3 v_2 - u_3 v_4 - u_4 v_2 + u_4 v_3) + x_2 (v_1 w_3 - v_1 w_4 - v_3 w_1 + v_3 w_4 + v_4 w_1 - v_4 w_3) \\
    &+ y_2 (-u_1 w_3 + u_1 w_4 + u_3 w_1 - u_3 w_4 - u_4 w_1 + u_4 w_3) + z_2 (u_1 v_3 - u_1 v_4 - u_3 v_1 + u_3 v_4 + u_4 v_1 - u_4 v_3) \\
    &+ x_3 (-v_1 w_2 + v_1 w_4 + v_2 w_1 - v_2 w_4 - v_4 w_1 + v_4 w_2) + y_3 (u_1 w_2 - u_1 w_4 - u_2 w_1 + u_2 w_4 + u_4 w_1 - u_4 w_2) \\
    &+ z_3 (-u_1 v_2 + u_1 v_4 + u_2 v_1 - u_2 v_4 - u_4 v_1 + u_4 v_2) + x_4 (v_1 w_2 - v_1 w_3 - v_2 w_1 + v_2 w_3 + v_3 w_1 - v_3 w_2) \\
    &+ y_4 (-u_1 w_2 + u_1 w_3 + u_2 w_1 - u_2 w_3 - u_3 w_1 + u_3 w_2) + z_4 (u_1 v_2 - u_1 v_3 - u_2 v_1 + u_2 v_3 + u_3 v_1 - u_3 v_2) \\
    q = &\ (-y_2 z_3 + y_2 z_4 + y_3 z_2 - y_3 z_4 - y_4 z_2 + y_4 z_3)^2 + (x_2 z_3 - x_2 z_4 - x_3 z_2 + x_3 z_4 + x_4 z_2 - x_4 z_3)^2 \\
    &+ (-x_2 y_3 + x_2 y_4 + x_3 y_2 - x_3 y_4 - x_4 y_2 + x_4 y_3)^2 + (y_1 z_3 - y_1 z_4 - y_3 z_1 + y_3 z_4 + y_4 z_1 - y_4 z_3)^2 \\
    &+ (-x_1 z_3 + x_1 z_4 + x_3 z_1 - x_3 z_4 - x_4 z_1 + x_4 z_3)^2 + (x_1 y_3 - x_1 y_4 - x_3 y_1 + x_3 y_4 + x_4 y_1 - x_4 y_3)^2 \\
    &+ (-y_1 z_2 + y_1 z_4 + y_2 z_1 - y_2 z_4 - y_4 z_1 + y_4 z_2)^2 + (x_1 z_2 - x_1 z_4 - x_2 z_1 + x_2 z_4 + x_4 z_1 - x_4 z_2)^2 \\
    &+ (-x_1 y_2 + x_1 y_4 + x_2 y_1 - x_2 y_4 - x_4 y_1 + x_4 y_2)^2 + (y_1 z_2 - y_1 z_3 - y_2 z_1 + y_2 z_3 + y_3 z_1 - y_3 z_2)^2 \\
    &+ (-x_1 z_2 + x_1 z_3 + x_2 z_1 - x_2 z_3 - x_3 z_1 + x_3 z_2)^2 + (x_1 y_2 - x_1 y_3 - x_2 y_1 + x_2 y_3 + x_3 y_1 - x_3 y_2)^2
\end{aligned}
\end{equation}
and with the constant-volume constraint:
\begin{equation}
\begin{aligned}
    0 = &\ (-y_2 z_3 + y_2 z_4 + y_3 z_2 - y_3 z_4 - y_4 z_2 + y_4 z_3) u_1 + (x_2 z_3 - x_2 z_4 - x_3 z_2 + x_3 z_4 + x_4 z_2 - x_4 z_3) v_1 \\
    &+ (-x_2 y_3 + x_2 y_4 + x_3 y_2 - x_3 y_4 - x_4 y_2 + x_4 y_3) w_1 + (y_1 z_3 - y_1 z_4 - y_3 z_1 + y_3 z_4 + y_4 z_1 - y_4 z_3) u_2 \\
    &+ (-x_1 z_3 + x_1 z_4 + x_3 z_1 - x_3 z_4 - x_4 z_1 + x_4 z_3) v_2 + (x_1 y_3 - x_1 y_4 - x_3 y_1 + x_3 y_4 + x_4 y_1 - x_4 y_3) w_2 \\
    &+ (-y_1 z_2 + y_1 z_4 + y_2 z_1 - y_2 z_4 - y_4 z_1 + y_4 z_2) u_3 + (x_1 z_2 - x_1 z_4 - x_2 z_1 + x_2 z_4 + x_4 z_1 - x_4 z_2) v_3 \\
    &+ (-x_1 y_2 + x_1 y_4 + x_2 y_1 - x_2 y_4 - x_4 y_1 + x_4 y_2) w_3 + (y_1 z_2 - y_1 z_3 - y_2 z_1 + y_2 z_3 + y_3 z_1 - y_3 z_2) u_4 \\
    &+ (-x_1 z_2 + x_1 z_3 + x_2 z_1 - x_2 z_3 - x_3 z_1 + x_3 z_2) v_4 + (x_1 y_2 - x_1 y_3 - x_2 y_1 + x_2 y_3 + x_3 y_1 - x_3 y_2) w_4
\end{aligned}
\end{equation}

\subsection{Conserved quantity results}
{\bf Known conserved quantities} Human experts obtained 8 conserved quantities for 2D and 12 conserved quantities for 3D. This is quite impressive already, given the complexity of the required calculations. However, they were unsure whether there are undiscovered conserved quantities, and whether the conserved quantities they found are in the simplest form.

{\bf SID} We run SID to search for conserved quantities. For the 2D case, SID finds there are indeed 8 conserved quantities, agreeing with human expertise. However, the conserved quantities found by SID are simpler. For example, all the conserved quantities discovered by SID are first or second-order polynomials, while the human expert finds a 4th order polynomial, which is later found (by human experts) to be derived from second-order conserved quantities discovered by SID. For the 3D case, SID finds there are 14 conserved quantities, i.e., two more than human experts discovered. Below we present SID's results given various basis sets (polynomials up to order $n=\{1,2,3,4\}$).

\subsubsection{2D}

\begin{figure}[htbp]
\includegraphics[width=0.49\linewidth]{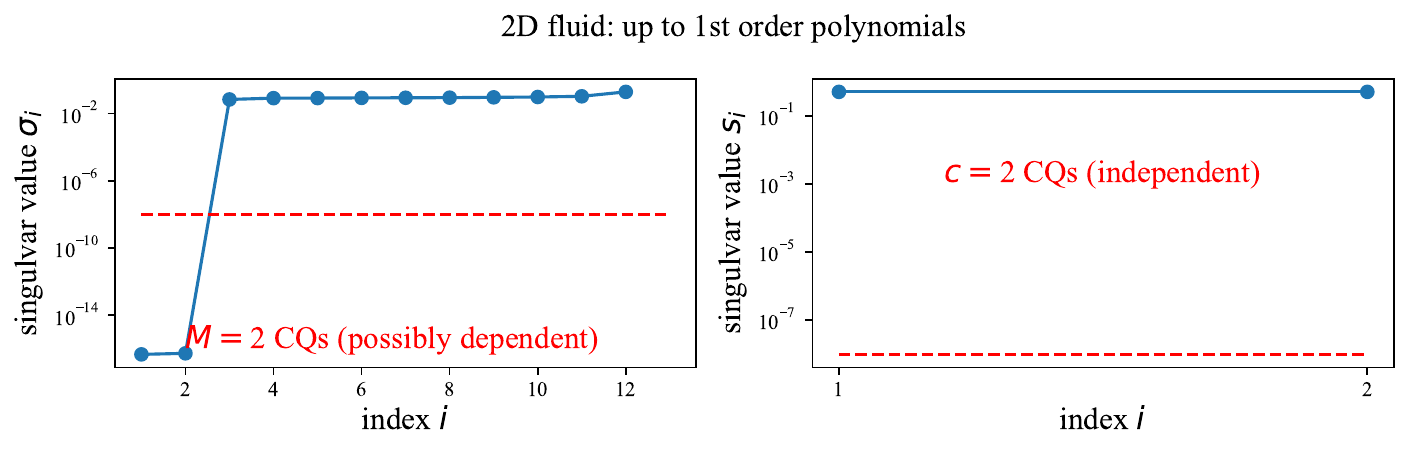} 
\includegraphics[width=0.49\linewidth]{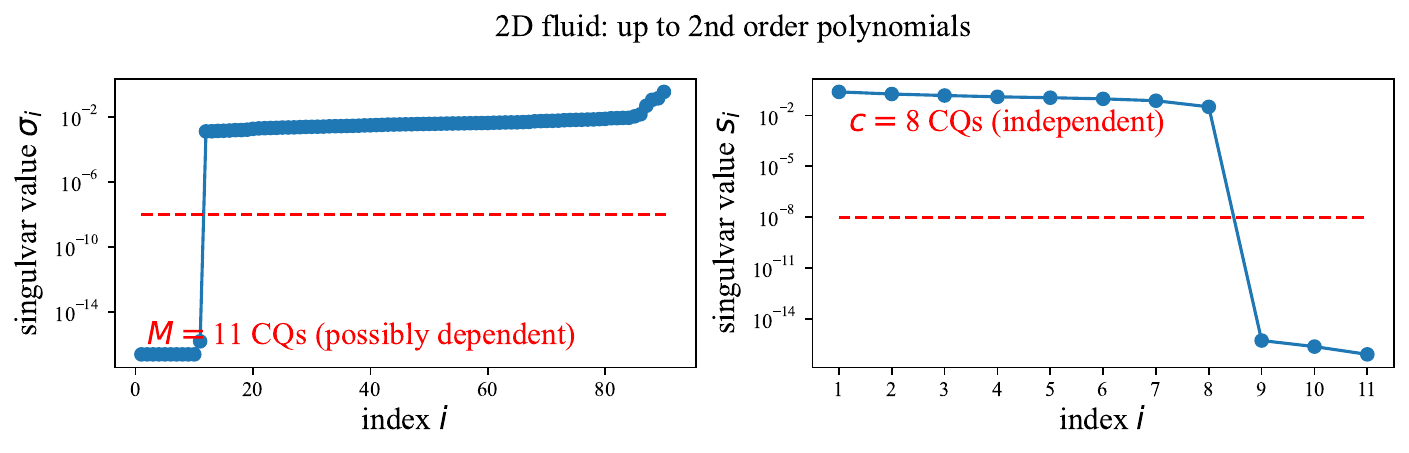} 
\includegraphics[width=0.49\linewidth]{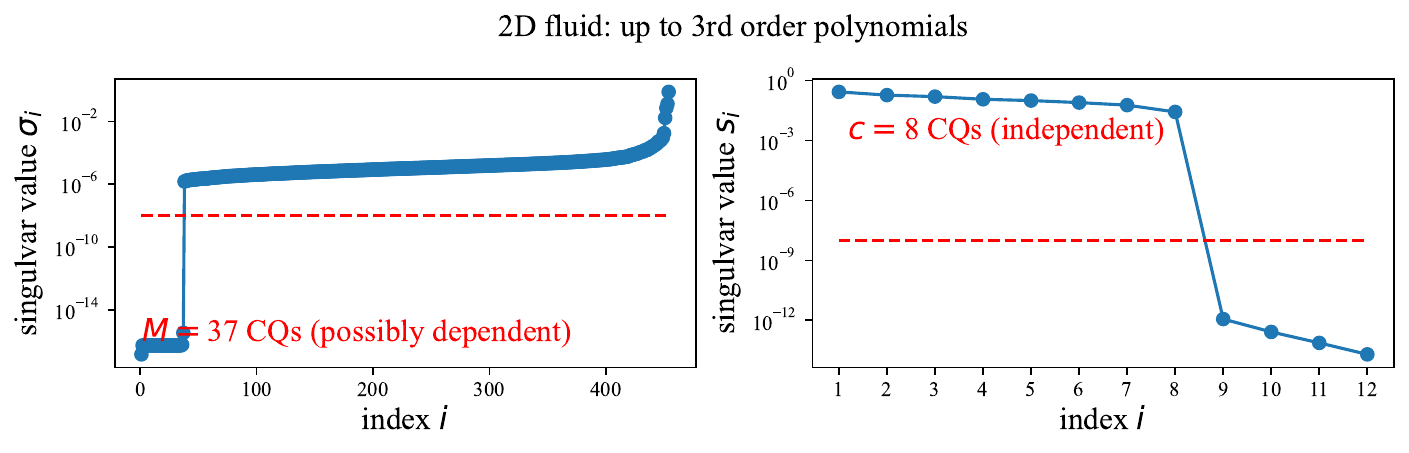}
\includegraphics[width=0.49\linewidth]{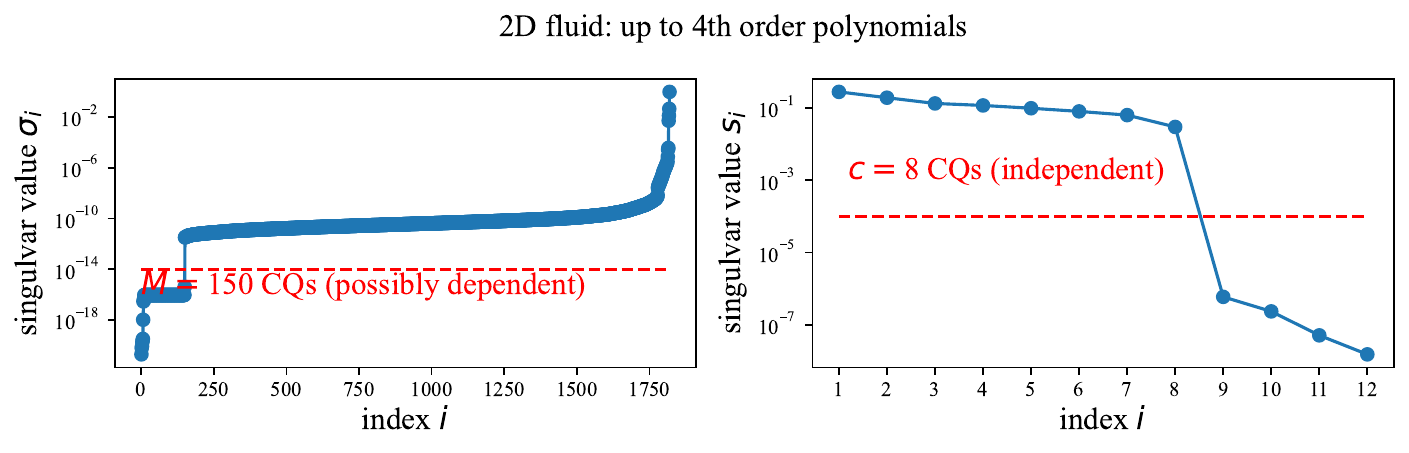} 
\caption{SID's results of the a simplified fluid system (2D) with various basis set (polynomials up to order $n$). When $n=1,2,3,4$, SID discovers 2, 8, 8, 8 conserved quantities respectively.}
\label{fig:fluid_2D_sv}
\end{figure}

\begin{figure}
    \centering
    \includegraphics[width=0.6\linewidth]{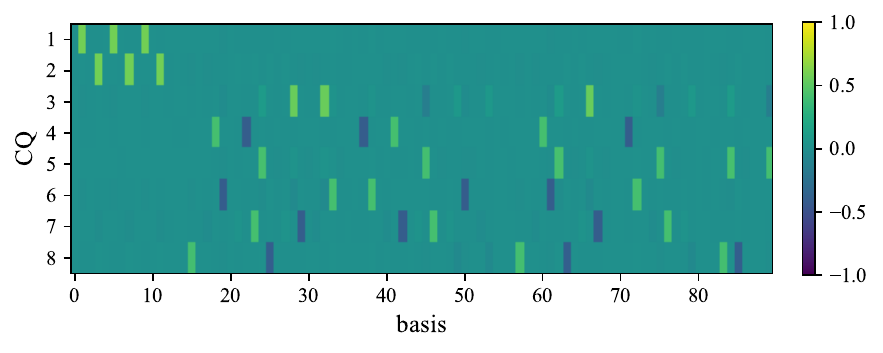}
    \caption{SID's discovered conserved quantities (coefficients) for the simplified fluid system (2D) with the basis set chosen to be polynomials up to order 2.}
    \label{fig:fluid_2D_CQ}
\end{figure}

As shown in Figure~\ref{fig:fluid_2D_sv}, when $n=1,2,3,4$, SID discovers 2, 8, 8, 8 conserved quantities, respectively. This means that there is no new 3rd order conserved quantity. Therefore we set $n=2$, and visualize the discovered conserved quantities (coefficients) in Figure~\ref{fig:fluid_2D_CQ}, which appear to be sparse (hence interpretable).

\subsubsection{3D}
\begin{figure}[htbp]
\includegraphics[width=0.49\linewidth]{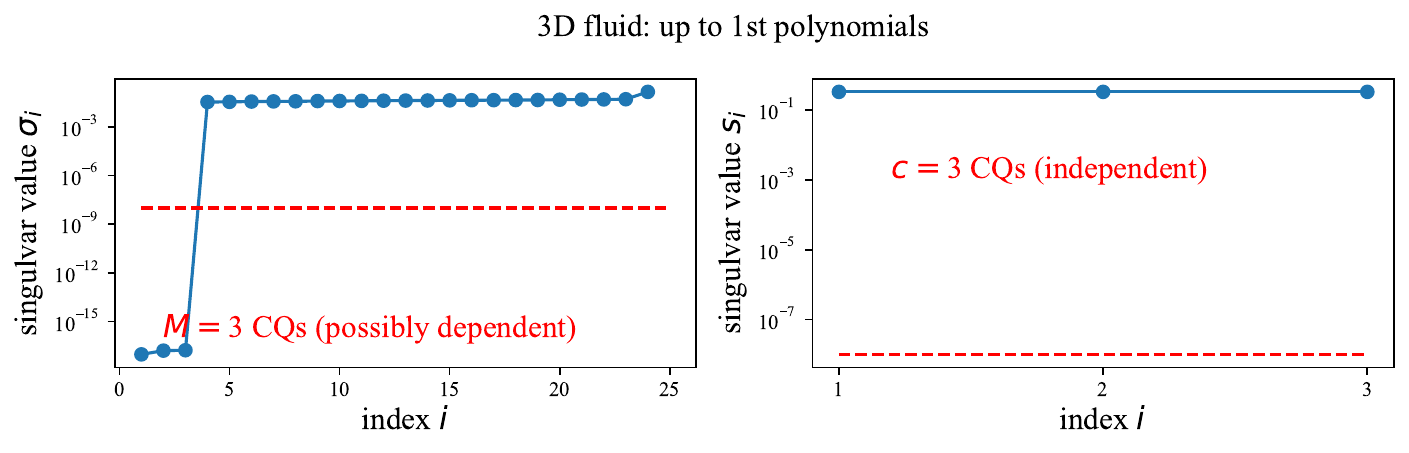} 
\includegraphics[width=0.49\linewidth]{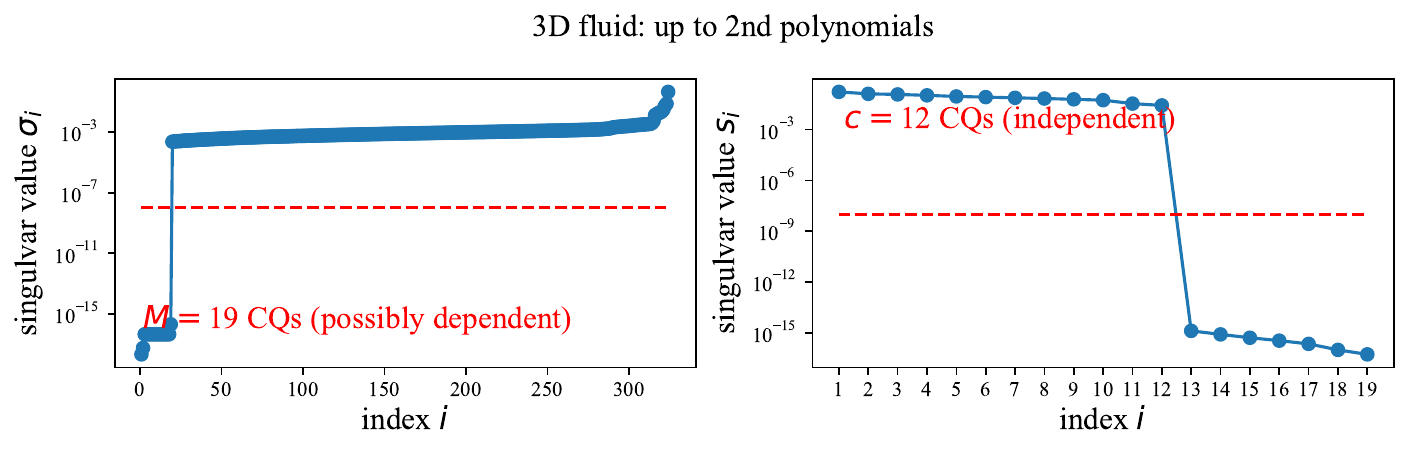} 
\includegraphics[width=0.49\linewidth]{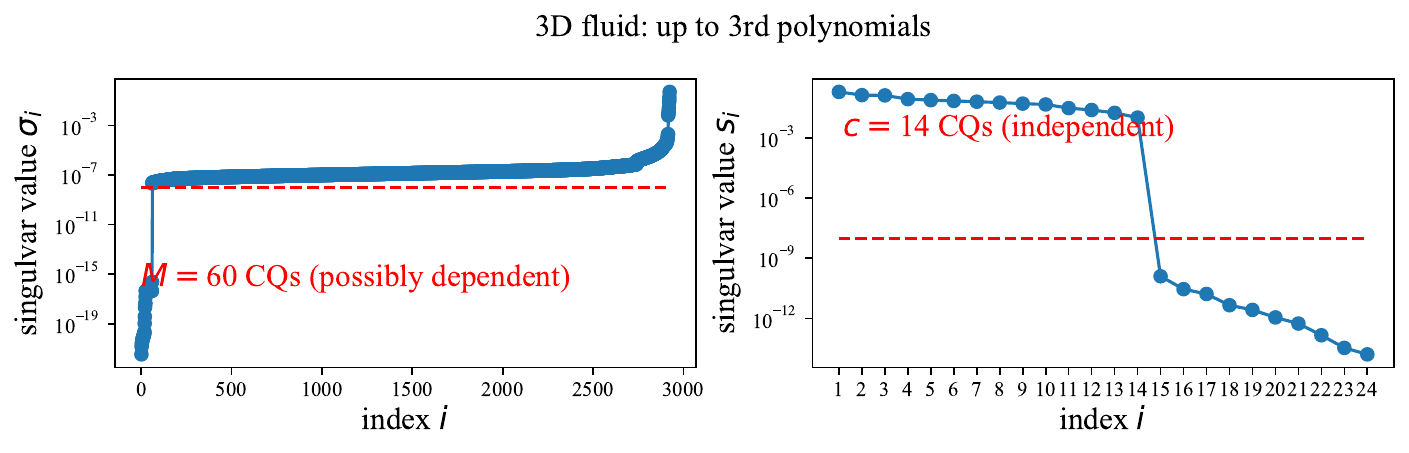}
\includegraphics[width=0.49\linewidth]{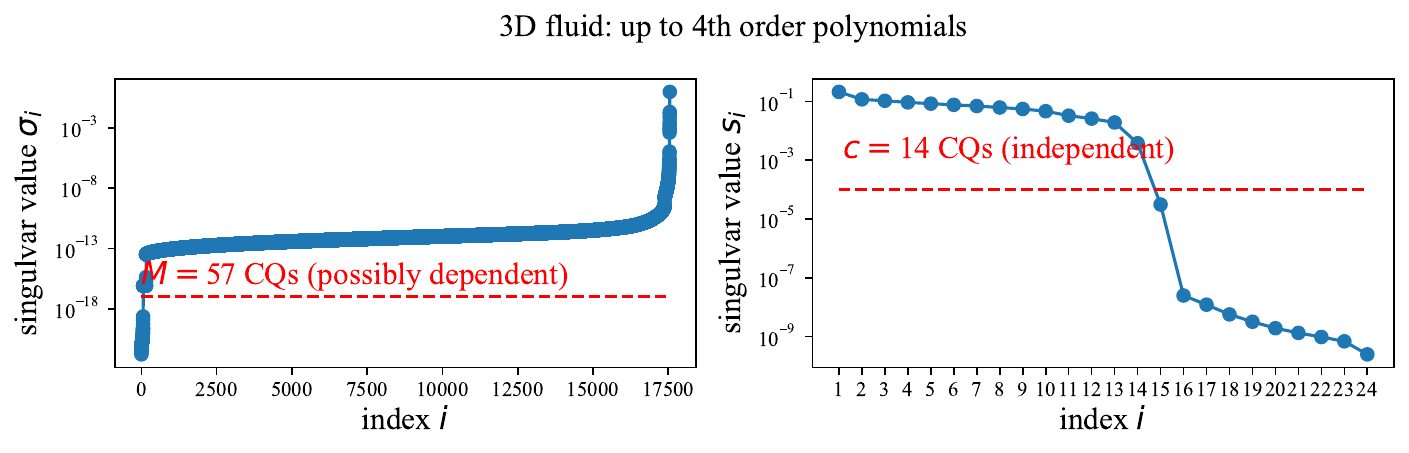}
\caption{SID's results of the a simplified fluid system (3D) with various basis set (polynomials up to order $n$). When $n=1,2,3,4$, SID discovers 3, 12, 14, 14 conserved quantities respectively.}
\label{fig:fluid_3D_sv}
\end{figure}

\begin{figure}[htbp]
    \centering
    \includegraphics[width=0.6\linewidth]{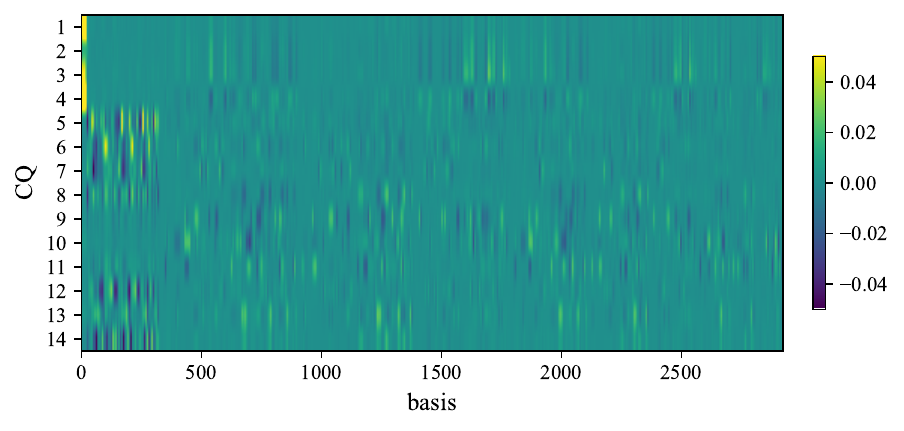}
    \caption{SID's discovered conserved quantities (coefficients) for the simplified fluid system (3D) with the basis set chosen to be polynomials up to order 3.}
    \label{fig:fluid_3D_CQ}
\end{figure}

As shown in Figure~\ref{fig:fluid_3D_sv}, when $n=1,2,3,4$, SID discovers 3, 12, 14, 14 conserved quantities, respectively. For $n=3$, we visualize SID's discovered conserved quantities (coefficients) in Figure~\ref{fig:fluid_3D_CQ}, which appear to be sparse (hence interpretable).

\subsection{Summary of conserved quantities (combining experts and SID)}
We have disclaimed frequently throughout the paper that SID by no means replaces human scientists, but rather is a helpful assistant that can facilitate the discovery process when an expert uses it. This is because: (1) SID needs the input of basis function set from human experts; (2) The physical interpretation of SID's results is still left to human experts. By combining the knowledge from human experts and SID's findings, we summarize the conserved quantities of the fluid system below.

\subsubsection{2D}

For the 2D case, the important quantities are:
\begin{equation}
\small
\begin{aligned}
& x_{cm}=\frac{1}{3}(x_1 + x_2 + x_3),\quad y_{cm}=\frac{1}{3}(y_1 + y_2 + y_3), \\
& \Bar{x}_1=x_1 - x_{cm},\quad \Bar{y}_1=y_1 - y_{cm}, \quad \Bar{x}_2=x_2 - x_{cm}, \quad \Bar{y}_2=y_2 - y_{cm}, \quad \Bar{x}_3=x_3 - x_{cm}, \quad \Bar{y}_3=y_3 - y_{cm}, \\
& u_{cm}=\frac{1}{3}(u_1 + u_2 + u_3), \quad v_{cm}=\frac{1}{3}(v_1 + v_2 + v_3), \\
& \Bar{u}_1=u_1 - u_{cm}, \quad \Bar{v}_1=v_1 - v_{cm},\quad \Bar{u}_2=u_2 - u_{cm},\quad \Bar{v}_2=v_2 - v_{cm},\quad \Bar{u}_3=u_3 - u_{cm},\quad \Bar{v}_3=v_3 - v_{cm}, \\
& L_{cm}=x_{cm}v_{cm} - y_{cm}u_{cm} \\
& L=\Bar{v}_1\Bar{x}_1 + \Bar{v}_2\Bar{x}_2 + \Bar{v}_3\Bar{x}_3 - \Bar{u}_1\Bar{y}_1 - \Bar{u}_2\Bar{y}_2 -  \Bar{u}_3\Bar{y}_3, \\
& E=\Bar{u}_1^2 + \Bar{v}_1^2 + \Bar{u}_2^2 + \Bar{v}_2^2 + \Bar{u}_3^2 + \Bar{v}_3^2, \\
& I=\Bar{x}_1^2 + \Bar{y}_1^2 + \Bar{x}_2^2 + \Bar{y}_2^2 + \Bar{x}_3^2 + \Bar{y}_3^2, \\
&J=(x_1-x_2)^2+(y_1-y_2)^2+(x_2-x_3)^2+(y_2-y_3)^2+(x_3-x_1)^2+(y_3-y_1)^2 , \\
& K=u_1 v_2 - u_1 v_3 - u_2 v_1 + u_2 v_3 + u_3 v_1 - u_3 v_2, \\
& A=x_1 y_2 - x_1 y_3 - x_2 y_1 + x_2 y_3 + x_3 y_1 - x_3 y_2, \\
& D=u_1(y_2 - y_3) + v_1(x_3 - x_2) + u_2(y_3 - y_1) + v_2(x_1 - x_3) + u_3(y_1 - y_2) + v_3(x_2 - x_1), \\
& \omega = u_1(x_2 - x_3) + u_2(x_3 - x_1) + u_3(x_1 - x_2) + v_1(y_2 - y_3)  + v_2(y_3 - y_1) + v_3(y_1 - y_2), \\
& G=\Bar{u}_1\Bar{x}_1 + \Bar{u}_2\Bar{x}_2 + \Bar{u}_3\Bar{x}_3 + \Bar{v}_1\Bar{y}_1 + \Bar{v}_2\Bar{y}_2 + \Bar{v}_3\Bar{y}_3
\end{aligned}
\end{equation}
The conserved quantities for the 2D case are
$$u_{cm}, v_{cm}, L_{cm}, L, E, A, D, \omega, IK$$ 
However, there is a relationship between these conserved quantities
$$IK = EA-L\omega-DG$$
which means there are only 8 independent conserved quantities. 

\begin{table}[ht]
\centering
\caption{Physical interpretations of important quantities (2D)}
\begin{tabular}[t]{ll}
\hline
Notation & Meaning\\
\hline
$u_{cm}$ & velocity of center of mass in x direction\\
$v_{cm}$ & velocity of center of mass in y direction\\
$L_{cm}$ & Angular momentum of center of mass\\
$L$ & Angular momentum with respect to center of mass\\
$E$ & Energy with respect to center of mass\\
$A$ & Area\\
$D$ & Divergence\\
$\omega$ & Vorticity\\\hline
$I$ & Moment of inertia with respect to center of mass\\
$G$ & Derivative of moment of inertia\\
\hline
\end{tabular}
\end{table}%

Global conserved quantities can be expressed in terms of integrals, where subscript $g$ stands for global and $\rho$ is density (suppose the support of the fluid is $\Omega\in\mathbb{R}^2$, the boundary is $\partial\Omega$):
\begin{equation}
    \begin{aligned}
        & \mat{p}_g= \iint_\Omega \rho\textbf{v} \,dA,\quad\textbf{v}\equiv (u,v),\quad {\rm (momentum, 2D\ vector)} \\
        & A_g = \iint_\Omega \,dA\quad {\rm (area, scalar)} \\
        & L_{cm,g} = \mat{r}_{cm}\times \mat{p}_g, \quad \mat{r}_{cm}\equiv \frac{1}{A_g}\iint_\Omega\mat{r}dA\quad  {\rm (momentum\ of\  COM, scalar)} \\
        & L_g = \left(\iint_\Omega \rho\textbf{r}\times\textbf{v} \,dA\right) - L_{cm,g}\quad {\rm (momentum\ in\  COM, scalar)}\\
        & E_g = \left(\iint_\Omega \rho|\textbf{v}|^2 \,dA\right)-\frac{|\mat{p}_g|^2}{2\rho A_g} \quad {\rm (energy\ in\ COM, scalar)} \\
        & D_g = \oint_{\partial\Omega} \,\mat{v}\cdot d\mat{n} \equiv 0 \quad {\rm (outflow, scalar)} \\
        & C_g = \iint_\Omega \omega dA = \oint_{\partial\Omega} \textbf{v} \cdot \,d \textbf{l}\quad {\rm (circulation, scalar)} \\
    \end{aligned}
\end{equation}

\subsubsection{3D}
For the 3D case, the important quantities are:

\begingroup
\allowdisplaybreaks
\begin{equation}
\small
\begin{aligned}
    & x_{cm}=\frac{1}{4}(x_1 + x_2 + x_3 +x_4),\quad y_{cm}=\frac{1}{4}(y_1 + y_2 + y_3 +y_4), \quad z_{cm}=\frac{1}{4}(z_1 + z_2 + z_3 +z_4) \\
    & \Bar{x}_1=x_1 - x_{cm},\quad \Bar{y}_1=y_1 - y_{cm}, \quad \Bar{z}_1=z_1 - z_{cm}, \\
    &\Bar{x}_2=x_2 - x_{cm},\quad \Bar{y}_2=y_2 - y_{cm},\quad \Bar{z}_2=z_2 - z_{cm}, \\
    & \Bar{x}_3=x_3 - x_{cm},\quad \Bar{y}_3=y_3 - y_{cm},\quad \Bar{z}_3=z_3 - z_{cm}, \\
    & \Bar{x}_4=x_4 - x_{cm}, \quad \Bar{y}_4=y_4 - y_{cm} ,\quad \Bar{z}_4=z_4 - z_{cm}, \\
    & u_{cm}=\frac{1}{4}(u_1 + u_2 + u_3 + u_4),\quad v_{cm}=\frac{1}{4}(v_1 + v_2 + v_3 + v_4), \quad w_{cm}=\frac{1}{4}(w_1 + w_2 + w_3 + w_4), \\
    & \Bar{u}_1=u_1 - u_{cm},\quad \Bar{v}_1=v_1 - v_{cm},\quad \Bar{w}_1=w_1 - w_{cm} \\
    & \Bar{u}_2=u_2 - u_{cm},\quad \Bar{v}_2=v_2 - v_{cm},\quad \Bar{w}_2=w_2 - w_{cm}, \\
    & \Bar{u}_3=u_3 - u_{cm},\quad \Bar{v}_3=v_3 - v_{cm},\quad \Bar{w}_3=w_3 - w_{cm}, \\
    & \Bar{u}_4=u_4 - u_{cm}, \quad \Bar{v}_4=v_4 - v_{cm},\quad \Bar{w}_4=w_4 - w_{cm}, \\
\end{aligned}
\end{equation}

\begin{equation}
\small
\begin{aligned}
    & A^x_1=-y_2 z_3 + y_2 z_4 + y_3 z_2 - y_3 z_4 - y_4 z_2 + y_4 z_3,\quad A^y_1=x_2 z_3 - x_2 z_4 - x_3 z_2 + x_3 z_4 + x_4 z_2 - x_4 z_3, \\
    & A^z_1=-x_2 y_3 + x_2 y_4 + x_3 y_2 - x_3 y_4 - x_4 y_2 + x_4 y_3,\quad A^x_2=y_1 z_3 - y_1 z_4 - y_3 z_1 + y_3 z_4 + y_4 z_1 - y_4 z_3, \\
    & A^y_2=-x_1 z_3 + x_1 z_4 + x_3 z_1 - x_3 z_4 - x_4 z_1 + x_4 z_3,\quad A^z_2=x_1 y_3 - x_1 y_4 - x_3 y_1 + x_3 y_4 + x_4 y_1 - x_4 y_3, \\
    & A^x_3=-y_1 z_2 + y_1 z_4 + y_2 z_1 - y_2 z_4 - y_4 z_1 + y_4 z_2,\quad A^y_3=x_1 z_2 - x_1 z_4 - x_2 z_1 + x_2 z_4 + x_4 z_1 - x_4 z_2 \\
    & A^z_3=-x_1 y_2 + x_1 y_4 + x_2 y_1 - x_2 y_4 - x_4 y_1 + x_4 y_2, \quad A^x_4=y_1 z_2 - y_1 z_3 - y_2 z_1 + y_2 z_3 + y_3 z_1 - y_3 z_2 \\
    & A^y_4=-x_1 z_2 + x_1 z_3 + x_2 z_1 - x_2 z_3 - x_3 z_1 + x_3 z_2,\quad A^z_4=x_1 y_2 - x_1 y_3 - x_2 y_1 + x_2 y_3 + x_3 y_1 - x_3 y_2, \\
    & B^x_1=-v_2 w_3 + v_2 w_4 + v_3 w_2 - v_3 w_4 - v_4 w_2 + v_4 w_3, \quad B^y_1=u_2 w_3 - u_2 w_4 - u_3 w_2 + u_3 w_4 + u_4 w_2 - u_4 w_3 \\
    & B^z_1=-u_2 v_3 + u_2 v_4 + u_3 v_2 - u_3 v_4 - u_4 v_2 + u_4 v_3, \quad B^x_2=v_1 w_3 - v_1 w_4 - v_3 w_1 + v_3 w_4 + v_4 w_1 - v_4 w_3, \\
    & B^y_2=-u_1 w_3 + u_1 w_4 + u_3 w_1 - u_3 w_4 - u_4 w_1 + u_4 w_3,\quad B^z_2=u_1 v_3 - u_1 v_4 - u_3 v_1 + u_3 v_4 + u_4 v_1 - u_4 v_3, \\
    & B^x_3=-v_1 w_2 + v_1 w_4 + v_2 w_1 - v_2 w_4 - v_4 w_1 + v_4 w_2, \quad B^y_3=u_1 w_2 - u_1 w_4 - u_2 w_1 + u_2 w_4 + u_4 w_1 - u_4 w_2 \\
    & B^z_3=-u_1 v_2 + u_1 v_4 + u_2 v_1 - u_2 v_4 - u_4 v_1 + u_4 v_2,\quad B^x_4=v_1 w_2 - v_1 w_3 - v_2 w_1 + v_2 w_3 + v_3 w_1 - v_3 w_2 \\
    & B^y_4=-u_1 w_2 + u_1 w_3 + u_2 w_1 - u_2 w_3 - u_3 w_1 + u_3 w_2,\quad B^z_4=u_1 v_2 - u_1 v_3 - u_2 v_1 + u_2 v_3 + u_3 v_1 - u_3 v_2, \\
\end{aligned}
\end{equation}

\begin{equation}
\small
\begin{aligned}
    & E = \Bar{u}_1^2 + \Bar{v}_1^2 + \Bar{w}_1^2 + \Bar{u}_2^2 + \Bar{v}_2^2 + \Bar{w}_2^2 + \Bar{u}_3^2 + \Bar{v}_3^2 + \Bar{w}_3^2 + \Bar{u}_4^2 + \Bar{v}_4^2 + \Bar{w}_4^2, \\
    &L_x = \Bar{w}_1\Bar{y}_1 + \Bar{w}_2\Bar{y}_2 + \Bar{w}_3\Bar{y}_3 + \Bar{w}_4\Bar{y}_4 - \Bar{v}_1\Bar{z}_1 - \Bar{v}_2\Bar{z}_2 - \Bar{v}_3\Bar{z}_3 - \Bar{v}_4\Bar{z}_4, \\
    &L_y = \Bar{u}_1\Bar{z}_1 + \Bar{u}_2\Bar{z}_2 + \Bar{u}_3\Bar{z}_3 + \Bar{u}_4\Bar{z}_4 - \Bar{w}_1\Bar{x}_1 - \Bar{w}_2\Bar{x}_2 - \Bar{w}_3\Bar{x}_3 - \Bar{w}_4\Bar{x}_4, \\
    & L_z = \Bar{v}_1\Bar{x}_1 + \Bar{v}_2\Bar{x}_2 + \Bar{v}_3\Bar{x}_3 + \Bar{v}_4\Bar{x}_4 - \Bar{u}_1\Bar{y}_1 - \Bar{u}_2\Bar{y}_2 - \Bar{u}_3\Bar{y}_3 - \Bar{u}_4\Bar{y}_4, \\
    & L^{cm}_x=w_{cm}y_{cm} - v_{cm}z_{cm}, \quad L^{cm}_y=u_{cm}z_{cm} - w_{cm}x_{cm}, \quad L^{cm}_z=v_{cm}x_{cm} - u_{cm}y_{cm}, \\
    & V = \frac{1}{3}(x_1 A^x_1 + y_1 A^y_1 + z_1 A^z_1 + x_2 A^x_2 + y_2 A^y_2 + z_2 A^z_2 + x_3 A^x_3 + y_3 A^y_3 + z_3 A^z_3 + x_4 A^x_4 + y_4 A^y_4 + z_4 A^z_4), \\
    & \omega_x = v_1 A^z_1 - w_1 A^y_1  + v_2 A^z_2 - w_2 A^y_1 + v_3 A^z_3 - w_3 A^y_3 + v_4 A^z_4 - w_4 A^y_4, \\
    & \omega_y = w_1 A^x_1 - u_1 A^z_1 + w_2 A^x_2 - u_2 A^z_1 + w_3 A^x_3 - u_3 A^z_3 + w_4 A^x_4 - u_4 A^z_4, \\
    & \omega_z = u_1 A^y_1 - v_1 A^x_1 + u_2 A^y_2 - v_2 A^x_2 + u_3 A^y_3 - v_3 A^x_3 + u_4 A^y_4 - v_4 A^x_4, \\
    & D = u_1 A^x_1 + v_1 A^y_1 + w_1 A^z_1 + u_2 A^x_2 + v_2 A^y_2 + w_2 A^z_2 + u_3 A^x_3 + v_3 A^y_3 + w_3 A^z_3 + u_4 A^x_4 + v_4 A^y_4 + w_4 A^z_4, \\
    & I=\Bar{x}_1^2 + \Bar{y}_1^2 + \Bar{z}_1^2 + \Bar{x}_2^2 + \Bar{y}_2^2  + \Bar{z}_2^2 + \Bar{x}_3^2 + \Bar{y}_3^2 + \Bar{z}_3^2 + \Bar{x}_4^2 + \Bar{y}_4^2 + \Bar{z}_4^2, \\
    & I_{xx}=\Bar{x}_1^2 + \Bar{x}_2^2 + \Bar{x}_3^2 + \Bar{x}_4^2,\quad I_{yy}=\Bar{y}_1^2 + \Bar{y}_2^2 + \Bar{y}_3^2 + \Bar{y}_4^2,\quad I_{zz}=\Bar{z}_1^2 + \Bar{z}_2^2 + \Bar{z}_3^2 + \Bar{z}_4^2, \\
    & I_{xy} = \Bar{x}_1\Bar{y}_1 + \Bar{x}_2\Bar{y}_2 + \Bar{x}_3\Bar{y}_3 + \Bar{x}_4\Bar{y}_4, \quad I_{yz} = \Bar{y}_1\Bar{z}_1 + \Bar{y}_2\Bar{z}_2 + \Bar{y}_3\Bar{z}_3 + \Bar{y}_4\Bar{z}_4,\quad I_{zx} = \Bar{z}_1\Bar{x}_1 + \Bar{z}_2\Bar{x}_2 + \Bar{z}_3\Bar{x}_3 + \Bar{z}_4\Bar{x}_4, \\
    & G = \Bar{u}_1\Bar{x}_1 + \Bar{v}_1\Bar{y}_1 + \Bar{w}_1\Bar{z}_1 + \Bar{u}_2\Bar{x}_2 + \Bar{v}_2\Bar{y}_2 + \Bar{w}_2\Bar{z}_2 + \Bar{u}_3\Bar{x}_3 + \Bar{v}_3\Bar{y}_3 + \Bar{w}_3\Bar{z}_3 + \Bar{u}_4\Bar{x}_4 + \Bar{v}_4\Bar{y}_4 + \Bar{w}_4\Bar{z}_4, \\
\end{aligned}
\end{equation}

\begin{equation}\small
\begin{aligned}
    & K=x_1 B^x_1 + y_1 B^y_1 + z_1 B^z_1 + x_2 B^x_2 + y_2 B^y_2 + z_2 B^z_2 + x_3 B^x_3 + y_3 B^y_3 + z_3 B^z_3 + x_4 B^x_4 + y_4 B^y_4 + z_4 B^z_4, \\
    & K_{xx} = x_1 B^x_1 + x_2 B^x_2 + x_3 B^x_3 + x_4 B^x_4,\quad K_{yy} = y_1 B^y_1 + y_2 B^y_2 + y_3 B^y_3 + y_4 B^y_4, \quad  K_{zz} = z_1 B^z_1 + z_2 B^z_2 + z_3 B^z_3 + z_4 B^z_4, \\
    & K_{xy} = x_1 B^y_1 + x_2 B^y_2 + x_3 B^y_3 + x_4 B^y_4, \quad K_{yx} = y_1 B^x_1 + y_2 B^x_2 + y_3 B^x_3 + y_4 B^x_4, \\
    & K_{yz} = y_1 B^z_1 + y_2 B^z_2 + y_3 B^z_3 + y_4 B^z_4,\quad K_{zy} = y_1 B^z_1 + y_2 B^z_2 + y_3 B^z_3 + y_4 B^z_4, \\
    & K_{zx} = z_1 B^x_1 + z_2 B^x_2 + z_3 B^x_3 + z_4 B^x_4,\quad K_{xz} = x_1 B^z_1 + x_2 B^z_2 + x_3 B^z_3 + x_4 B^z_4, \\
    & J=(A^x_1)^2 + (A^y_1)^2 + (A^z_1)^2 + (A^x_2)^2 + (A^y_2)^2 + (A^z_2)^2 + (A^x_3)^2 + (A^y_3)^2 + (A^z_3)^2 + (A^x_4)^2 + (A^y_4)^2 + (A^z_4)^2, \\
    &C_1 = u_2 x_3 - u_2 x_4 - u_3 x_2 + u_3 x_4 + u_4 x_2 - u_4 x_3 + v_2 y_3 - v_2 y_4 - v_3 y_2 + v_3 y_4 + v_4 y_2 - v_4 y_3 + w_2 z_3 - w_2 z_4 - w_3 z_2 + w_3 z_4 + w_4 z_2 - w_4 z_3, \\
    &C_2 = -u_1 x_3 + u_1 x_4 + u_3 x_1 - u_3 x_4 - u_4 x_1 + u_4 x_3 - v_1 y_3 + v_1 y_4 + v_3 y_1 - v_3 y_4 - v_4 y_1 + v_4 y_3 - w_1 z_3 + w_1 z_4 + w_3 z_1 - w_3 z_4 - w_4 z_1 + w_4 z_3, \\
    & C_3 = u_1 x_2 - u_1 x_4 - u_2 x_1 + u_2 x_4 + u_4 x_1 - u_4 x_2 + v_1 y_2 - v_1 y_4 - v_2 y_1 + v_2 y_4 + v_4 y_1 - v_4 y_2 + w_1 z_2 - w_1 z_4 - w_2 z_1 + w_2 z_4 + w_4 z_1 - w_4 z_2, \\
    & C_4 = -u_1 x_2 + u_1 x_3 + u_2 x_1 - u_2 x_3 - u_3 x_1 + u_3 x_2 - v_1 y_2 + v_1 y_3 + v_2 y_1 - v_2 y_3 - v_3 y_1 + v_3 y_2 - w_1 z_2 + w_1 z_3 + w_2 z_1 - w_2 z_3 - w_3 z_1 + w_3 z_2, \\
\end{aligned}
\end{equation}
\endgroup

The conserved quantities of the 3D case are
$$u_{cm}, v_{cm}, w_{cm}, L^{cm}_x, L^{cm}_y,L^{cm}_z, L_x, L_y, L_z, E, V, D, C_1, C_2, C_3, C_4$$
The relationships between the conserved quantities are
$$C_1 + C_2 + C_3 + C_4=0, \quad u_{cm} L^{cm}_x + v_{cm} L^{cm}_y + w_{cm} L^{cm}_z=0,$$
which means there are only 14 independent conserved quantities.

Global conserved quantities can be expressed in terms of integrals, where subscript $g$ stands for global and $\rho$ is density (suppose the support of the fluid is $\Omega\in\mathbb{R}^3$, the boundary is $\partial\Omega$, $C$ stands for any closed loop on $\partial\Omega$):
\begin{equation}
    \begin{aligned}
        & V_g = \iiint_\Omega \,dV\quad {\rm (volume, scalar)} \\
        & \mat{p}_g= \iiint_\Omega \rho\textbf{v} \,dV,\quad  \textbf{v}=(u,v,w),\quad {\rm (momentum, 3D\ vector)} \\
        & E_g = \left(\iiint_\Omega \rho|\textbf{v}|^2 \,dV\right)-\frac{|\mat{p}_g|^2}{2\rho V_g} \quad {\rm (energy\ in\ COM, scalar)} \\
        & L_{cm,g} = \mat{r}_{cm}\times \mat{p}_g \quad \mat{r}_{cm}\equiv \frac{1}{V_g}\iiint_\Omega\mat{r}dV\quad  {\rm (momentum\ of\  COM, 3D\  vector)} \\
        & L_g = \left(\iiint_\Omega \rho\textbf{r}\times\textbf{v} \,dV\right) - L_{cm,g}\quad {\rm (momentum\ in\  COM, 3D\  vector)}\\
        & D_g = \oiint_{\partial\Omega} \,\mat{v}\cdot d\mat{n} \equiv 0 \quad {\rm (outflow, scalar)} \\
        & C_g = \oint_{C\in\partial \Omega} \textbf{v} \cdot \,d \textbf{l}\quad {\rm (circulation, scalar)} \\
    \end{aligned}
\end{equation}

\begin{table}[ht]
\centering
\caption{Physical interpretations of important quantities (3D)}
\begin{tabular}[t]{ll}
\hline
Notation & Meaning\\
\hline
$u_{cm}$ & velocity of center of mass in x direction\\
$v_{cm}$ & velocity of center of mass in y direction\\
$w_{cm}$ & velocity of center of mass in z direction\\
$L^{cm}_x$ & Angular momentum of center of mass in x direction\\
$L^{cm}_y$ & Angular momentum of center of mass in y direction\\
$L^{cm}_z$ & Angular momentum of center of mass in z direction\\
$L_x$ & Angular momentum with respect to center of mass in x direction\\
$L_y$ & Angular momentum with respect to center of mass in y direction\\
$L_z$ & Angular momentum with respect to center of mass in z direction\\
$E$ & Energy with respect to center of mass\\
$V$ & Volume\\
$D$ & Divergence\\
$C_i$ & Circulations $i=1,2,3,4$\\\hline
$\mat{I}$ & Moment of inertia with respect to center of mass\\
$\mat{\omega}$ & Vorticity\\
$G$ & Derivative of moment of inertia\\
\hline
\end{tabular}
\end{table}%

\end{document}